\documentclass{article}
\usepackage[T1]{fontenc}
\usepackage{graphicx,amsmath,amsfonts,amssymb,color,amsthm,psfrag}
\usepackage{here}

\newtheorem{thm}{Theorem}[section]

\newtheorem{lem}[thm]{Lemma}

\newtheorem{prop}[thm]{Proposition}

\newtheorem{rem}[thm]{Remark}

\numberwithin{equation}{section}

\renewcommand{\d}[1]{\, d  #1}

\renewcommand{\theta}{\vartheta}

\title{Non uniqueness of stationary measures for self-stabilizing processes}
\author{S. Herrmann and  J. Tugaut\\
  Institut de Math\'ematiques Elie Cartan - UMR 7502\\
  Nancy-Universit\'e, CNRS, INRIA\\
  B.P. 239, 54506 Vandoeuvre-l\`es-Nancy Cedex, France}
\begin{document}
\newcommand{\uem}{U_{\epsilon,\mu}}
\newcommand{\xem}{x_{\epsilon,\mu}}
\maketitle
\begin{abstract}
We investigate the existence of invariant measures for self-stabilizing diffusions. These stochastic processes represent roughly the behavior of some Brownian particle moving in a double-well landscape and attracted by its own law. This specific self-interaction leads to nonlinear stochastic differential equations and permits to point out singular phenomenons like non uniqueness of associated stationary measures. The existence of several invariant measures is essentially based on the non convex environment and requires generalized Laplace's method approximations. 
\end{abstract}
\medskip

{\bf Key words and phrases:} self-interacting diffusion;
stationary measures; double well potential; perturbed
dynamical system; Laplace's method; fixed point theorem.\par\medskip

{\bf 2000 AMS subject classifications:} primary 60H10;
secondary: 60J60, 60G10, 41A60\par\medskip

\section{Introduction}
The aim of this paper is to present some new and surprising results concerning the existence of invariant probability measures for one-dimensional self-stabilizing diffusions. The specificity of such diffusion is the attraction of its paths by the own law of the stochastic process. The dynamical system solved by self-stabilizing diffusions can be characterized by three essential elements: first the system is governed by a double-well potential $V$ which represents roughly the environment of the process, secondly some interaction potential $F$ describes how strong the attraction between the process and its own law is, and finally the system is perturbed by some Brownian motion with small amplitude $(\sqrt{\epsilon} B_t,\ t\ge 0)$. \\
Let us denote by $u_t^{\epsilon}(dx)$ the law of the self-stabilizing diffusion $(X_t^\epsilon,\ t\ge 0)$, then the SDE satisfied by $(X_t^\epsilon)$ is given by:
\[
X_t^{\epsilon}=X_0+\sqrt{\epsilon}B_t-\int_0^tV'(X_s^{\epsilon})ds-\int_0^t\int_{\mathbb{Re}}F'(X_s^{\epsilon}-x)du_s^{\epsilon}(x)ds, \quad \epsilon>0 . \quad\quad\hfill{(E^{\epsilon,X_0})}
\]
Introducing the notation of the convolution product, 
$(E^{\epsilon,X_0})$ can be written as follows:
\begin{equation}
\label{eq:equation_depart}
X_t^{\epsilon}=X_0+\sqrt{\epsilon}B_t-\int_0^t\left(V'+F'\ast u_s^{\epsilon}\right)\left(X_s^{\epsilon}\right)ds.
\end{equation}
Let us just note that the interaction part of the drift term is related to the diffusion in some simple way: $F'\ast u^\epsilon_t(x)=\mathbb{E}[F'(x-X_t^\epsilon)]$. This way of characterizing the drift term essentially points out the structure of the attraction between the paths of the diffusion and its law. Self-interaction corresponds obviously to \emph{mean fields} stabilization.

Self-stabilizing diffusion paths can usually be approximated by the movement of some specific Brownian particle belonging to a huge ensemble of identical ones. In this global system each particle is submitted to the same forces. First it moves in the potential landscape characterized by the double-well function $V$ and accordingly it is attracted by positions which minimize the potential. The second force which acts on the system is the interaction between all the particles. More precisely each one is attracted by all the others. This attraction can for instance be thought of as being generated by electromagnetic effects. In this case, the solution of the global system doesn't represent some spatial position but some electromagnetic charge.

The huge particle system containing $N$ elements is governed by the following
stochastic differential equation
\begin{align} \notag
  dX_t^{i,N} &= \sqrt{\varepsilon} \,dW_t^i, -V'(X_t^{i,N}) \,dt
               - \frac{1}{N} \sum_{j=1}^N F'(X_t^{i,N} - X_t^{j,N}) \,dt,
               \\
  X_0^{i,N} &= x_0\in\mathbb{Re},\quad 1\le i\le N.
  \label{eq:particle_system}
\end{align}
Here the $W^i$ are independent Brownian motions. In the limit, as $N$ becomes large, the interaction part of the drift term is approximatively the average with respect to the law of one characteristic particle of the system (law of large numbers framework). More precisely the empirical measures $\frac{1}{N} \sum_{j=1}^N
\delta_{X_t^{j,N}}$ converges to some law $u_t^\varepsilon$
for each fixed time and noise intensity, and each individual
particle's motion converges in probability to the solution of the
diffusion equation
\begin{equation} \label{eq00}
  d X_t^i =  \sqrt{\varepsilon} d W_t^i-V'(X_t^i) \,dt - \int_{\mathbb{Re}} F'(X_t^i - x) \d{u_t^\varepsilon(x)} \,dt.
\end{equation}
Interacting particle systems such as~\eqref{eq:particle_system} have
been studied from various points of view. A survey about the general
setting for interaction (under global Lipschitz and boundedness
assumptions) may be found in~\cite{Sznitman}.

The aim of this paper is to consider both the existence and the uniqueness of stationary measures for the self-stabilizing diffusion ($E^{\epsilon,X_0}$). In \cite{herrmann08} Herrmann, Imkeller and Peithmann proved the existence of some unique strong solution to equation \eqref{eq:equation_depart} generalizing previous results obtained by Benachour, Roynette, Talay and Vallois \cite{BRTV} in the context of constant environment potential $V$ ($V'(x)=0$ for all $x\in\mathbb{Re}$). We choose their work as basis for developing our study. Nevertheless there exist several different papers dealing with the existence problem for self-stabilizing diffusion, each of them concerning other families of interaction functions. Let us cite 
McKean who studied in some earlier work a class of Markov processes that contains the
solution of the limiting equation under restrictive global Lipschitz assumptions
for the interaction~\cite{McKean}, Stroock and Varadhan who considered some local form of interaction \cite{Stroock}, Oelschl\"{a}ger who studied the particular case
where interaction is represented by the derivative of the Dirac
measure at zero \cite{Oelschl} and finally Funaki who addressed existence and
uniqueness for the martingale problem associated with
self-stabilizing diffusions~\cite{Funaki}.

Let us now focus our attention to the stationary measures. In \cite{BRTV} the authors emphasize that the invariant measure, corresponding to some given average, is unique in this particular constant potential $V$ situation. This feature is essential for further developments. The natural convergence question between the law of the process and the invariant measure, as time elapses, can then be analyzed, see \cite{BRV}. This kind of convergence was also considered by Tamura under different assumptions on the structure of the interaction, see \cite{Tamura1} and \cite{Tamura2}.

The presence of some potential gradient which describes the environment of the self-stabilizing diffusion is essential for the question of existence and uniqueness of invariant measures. In particular, if the landscape is represented by some symmetric double-well potential then surprising effects appear due to the lack of convexity: we shall prove that, under suitable conditions, there exist at least three invariant measures of which one is symmetric (Theorem \ref{thm:exist-sym}) and two are asymmetric or so-called \emph{outlying} (Theorem \ref{thm:exist:asym}). In the particular linear interaction case ($F'(x)=\alpha x$ with $\alpha>0$), these three measures constitute the whole set of invariant measures (Theorem \ref{prop:convex}) provided that $V''$ is a convex function.

The material of this paper is organized as follows: first we list several assumptions concerning both the interaction function $F$ and the environment potential $V$ which permit in particular to assure the existence of the self-stabilizing diffusion ($E^{\epsilon,X_0}$). In Section 2 preliminary results concerning the structure of the invariant measure (if it exists !) are developed. These results are essential for the construction of such measures. The question of existence starts to be addressed in Section 3 in the particular linear interaction context. After pointing out some symmetric and asymmetric invariant measures, we point out some nice context for which the whole set of stationary  measures can be described. This study is finally extended to the general interaction case in the last section.
We postpone different tools concerning asymptotic analysis based on Laplace's method to the Appendix.

\subsection{Main assumptions}
In order to study invariant measures for self-stabilizing diffusions, we especially need that \eqref{eq:equation_depart} admits some unique strong solution. For this reason, we assume that both the potential landscape $V$ and the interaction function $F$ satisfy some growth conditions and some regularity properties. Moreover we add some technical assumptions which permit to simplify the statements.\\
We assume the following properties for the function $V$:\\[10pt]
\begin{minipage}{7.5cm}
\begin{description}
  \item[(V-1)] \emph{Regularity:} $V\in \mathcal{C}^{\infty}(\mathbb{Re},\mathbb{Re})$. $\mathcal{C}^\infty$ denotes the Banach space of infinitely bounded continuously differentiable function.
  \item[(V-2)] \emph{Symmetry:} $V$ is an even function.
    \item[(V-3)] $V$ is a double-well potential. The equation $V'(x)=0$ admits exactly three solutions : $a$, $-a$ and $0$ with $a>0$;  $V''(a)>0$ and $V''(0)<0$. The bottoms of the wells are reached for $x=a$ and $x=-a$.
    \item[(V-4)] There exist two constants $C_4,C_2>0$ such that $\forall x\in\mathbb{Re}$, $V(x)\geq C_4x^4-C_2x^2$.
    \end{description}
    \end{minipage}
     \begin{minipage}{4.5cm}
     \begin{figure}[H]
     \psfrag{$V$}{$V$}
     \psfrag{$-a$}{$-a$}
     \psfrag{$a$}{$a$}
     \centerline{\includegraphics[width=3.5cm]{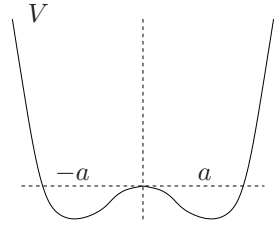}}
     \caption{Potential $V$}
     \end{figure}
     \end{minipage}
     \begin{description}
    \item[(V-5)] $\displaystyle\lim_{x\to\pm\infty}V''(x)=+\infty$ and $\forall x\geq a$, $V''(x)>0$.
    \item[(V-6)] The growth of the potential $V$ is at most polynomial: there exist $ q\in\mathbb{N}^*$ and $ C_q>0$ such that $\left|V'(x)\right|\leq C_q\left(1+x^{2q}\right)$.
    \item[(V-7)] \emph{Initialization:} $V(0)=0$.
\end{description}
Typically, $V$ is a double-well polynomial function. But our results can be applied to more general functions: regular functions with polynomial growth as $\vert x\vert$ becomes large.
We introduce the parameter $\theta$ which plays some important role in the following:
\begin{equation}
\label{eq:param_theta}
\theta=\sup_{x\in\mathbb{Re}}-V''(x).
\end{equation}
Let us note that the simplest example (most famous in the literature) is $V(x)=\frac{x^4}{4}-\frac{x^2}{2}$ which bottoms are localized in $-1$ and $1$ and with parameter $\theta=1$.\\
Let us now present the assumptions concerning the attraction function $F$.
\begin{description}
  \item[(F-1)] $F$ is an even polynomial function. Indeed we consider some classical situation: the attraction between two points $x$ and $y$ only depends on the distance $F(x-y)=F(y-x)$.
    \item[(F-2)] $F$ is a convex function.
    \item[(F-3)] $F'$ is a convex function on $\mathbb{Re}_+$ therefore for any $x\ge 0$ and $y\ge 0$ such that $x\geq y$ we get $F'(x)-F'(y)\geq F''(0)(x-y)$.
    \item[(F-4)] The polynomial growth of the attraction function $F$ is related to the growth condition (V-6): $|F'(x)-F'(y)|\leq C_q|x-y|(1+|x|^{2q-2}+|y|^{2q-2})$.
\end{description}
Let us define the parameter $\alpha\ge 0$ which shall play some essential role in following:
\begin{equation}\label{def:alpha}
F'(x)=\alpha x+F_0'(x)\quad\mbox{with}\ \alpha=F''(0)\geq0.
\end{equation}
In \cite{herrmann08}, Herrmann, Imkeller and Peithmann present sufficient conditions for the SDE \eqref{eq:equation_depart} to admit a unique strong solution. In particular, if $\mathbb{E}[X_0^{8q^2}]<+\infty$, with $q$ defined in (V-6), and if all the main assumptions just defined are satisfied, the existence and uniqueness of the solutions are proved. In the following we will always assume that the $(8q^2)$-th moment of the initial value $X_0$ is finite. This permits to study further the self-stabilizing diffusions and exhibit invariant measures.
\section{General structure of the invariant measures}
This section deals with different preliminary results describing the main structure of the invariant measures of $\left(E^{\epsilon,X_0}\right)$. First of all, there is some classical link between the stochastic differential equation and the associated parabolic partial differential equation which permits to characterize stationary measures.
\begin{lem}
\label{lem:PDE}
Let $u^{\epsilon}_t(x)$ denote the density of $(X_t^{\epsilon};\ t\ge 0)$ with respect to the Lebesgue measure. Then $u^\epsilon$ is solution of the following PDE:
\begin{eqnarray}
\label{eq:PDE}
\frac{\partial }{\partial t}u_t^\epsilon(x)=\frac{\epsilon}{2}\frac{\partial ^2}{\partial x^2}u_t^\epsilon(x)+\frac{\partial}{\partial x}\Big[u_t^\epsilon(x)\Big(V'(x)+\left(F'\ast u^\epsilon\right)(t,x)\Big)\Big]
\end{eqnarray}
for all $t>0$, $x\in\mathbb{Re}$ and $u_0^\epsilon(dx)=\mathbb{P}\left(X_0\in dx\right)$.\\ We recall that $\int_{\mathbb{R}}x^{8q^2}u_0^\epsilon(dx)<\infty$.
\end{lem}
\begin{proof}
Let $f\in\mathcal{C}^2(\mathbb{Re},\mathbb{Re})$ such that
\begin{eqnarray}\label{eq:cond_function}
\lim_{x\longrightarrow\pm\infty}f(x)=\lim_{x\longrightarrow\pm\infty}f'(x)=0.
\end{eqnarray}
By Itô's formula, we obtain
\begin{eqnarray*}
\mathbb{E}\left[f(X_t^{\epsilon})\right]&=&\mathbb{E}\left[f(X_0^{\epsilon})\right]+\mathbb{E}\left[\int_0^tf'(X_s^{\epsilon})\sqrt{\epsilon}dB_s\right]\\
&-&\mathbb{E}\left[\int_0^tf'(X_s^{\epsilon})\Big(V'(X_s^{\epsilon})+F'\ast
u^{\epsilon}_s(X_s^{\epsilon})\Big)ds+\frac{\epsilon}{2}\int_0^tf''(X_s^{\epsilon})ds\right].
\end{eqnarray*}
Taking the time derivative, we get
\begin{eqnarray*}\frac{d}{dt}\mathbb{E}\left[f(X_t^{\epsilon})\right]&=&-\mathbb{E}\left[f'(X_t^{\epsilon})\Big(V'(X_t^{\epsilon})+F'\ast u^{\epsilon}_t(X_t^{\epsilon})\Big)+\frac{\epsilon}{2}f''(X_t^{\epsilon})\right]\\
&=&-\int_{\mathbb{Re}}f'(x)\Big(V'(x)+F'\ast u^{\epsilon}_t(x)\Big)u_t^{\epsilon}(x)dx-\int_{\mathbb{Re}}\frac{\epsilon}{2}f''(x)u_t^{\epsilon}(x)dx.
\end{eqnarray*}
Since $f$ is a $\mathcal{C}^2$ function, integration by parts leads to
$$
\frac{d}{dt}\mathbb{E}\left[f(X_t^{\epsilon})\right]=\int_{\mathbb{Re}}f(x)\left\{\frac{\partial}{\partial x}\Big[\left(V'(x)+F'\ast u_t^{\epsilon}(x)\right)u_t^{\epsilon}(x)\Big]+\frac{\epsilon}{2}\frac{\partial^2}{\partial x^2}u_t^{\epsilon}(x)\right\}dx.
$$
Using the equality:
$$
\frac{d}{dt}\int_{\mathbb{Re}}f(x)u_t^{\epsilon}(x)dx=\int_{\mathbb{Re}}f(x)\frac{\partial}{\partial t}u_t^{\epsilon}(x)dx,
$$
we deduce for all $\mathcal{C}^2$-functions satisfying \eqref{eq:cond_function}:
$$\int_{\mathbb{Re}}f(x)\frac{\partial u_t^{\epsilon}(x)}{\partial t}\,dx=\int_{\mathbb{Re}}f(x)\left\{\frac{\partial}{\partial x}\Big[\left(V'(x)+F'\ast u_t^{\epsilon}(x)\right)u_t^{\epsilon}(x)\Big]+\frac{\epsilon}{2}\frac{\partial^2u_t^{\epsilon}(x)}{\partial x^2}\right\}dx$$
We obtain \eqref{eq:PDE} by identification.
\end{proof}
The density of $(X^\epsilon_t,\, t\ge 0)$ with respect to the
Lebesgue measure is solution to the parabolic PDE \eqref{eq:PDE}
(non-linear Kolmogorov equation): this implies in particular that
any stationary measure (if it exists !) satisfies some elliptic
 differential equation. This link between non-linear
differential equations and self-stabilizing diffusions permits to
express the invariant measure in some exponential form.
\begin{lem}
\label{lem:exp_form}
If there exists an invariant measure $u_{\epsilon}$ to $(E^{\epsilon,X_0})$ whose $(8q^2)$-moment is finite, then:
\begin{eqnarray}
\label{eq:lem:structure}
u_{\epsilon}(x)&=&\frac{1}{\lambda(u_{\epsilon})}\exp\left[-\frac{2}{\epsilon}\left(\int_0^xF'\ast u_{\epsilon}(y)dy+V(x)\right)\right]\\
&=&\frac{1}{\lambda(u_{\epsilon})}\exp\left[-\frac{2}{\epsilon}\Big(F\ast u_{\epsilon}(x)-F\ast u_{\epsilon}(0)+V(x)\Big)\right],\nonumber
\end{eqnarray}
where $\lambda(u_{\epsilon})$ denotes the normalization factor: $\int_\mathbb{Re} u_\epsilon(x)dx=1$.
\end{lem}
\begin{proof}
By \eqref{eq:PDE}, any stationary measure $u_{\epsilon}$ satisfies
\begin{eqnarray*}
\frac{\epsilon}{2}u_{\epsilon}''(x)+\Big(u_{\epsilon}(x)\left(V'(x)+F'\ast u_{\epsilon}(x)\right)\Big)'=0,\quad \mbox{for all}\ x\in\mathbb{Re}.
\end{eqnarray*}
By integrating the previous equality, we obtain the existence of some constant $C_{\epsilon}\in\mathbb{Re}$ such that
\begin{eqnarray*}
\frac{\epsilon}{2}u_{\epsilon}'(x)+u_{\epsilon}(x)(V'(x)+F'\ast u_{\epsilon}(x))=C_{\epsilon},\quad \mbox{for all}\ x\in\mathbb{Re}.
\end{eqnarray*}
Using the method of variation of parameters, the solution $u_{\epsilon}$ takes the following form
$$
u_{\epsilon}(x)=\Lambda_{\epsilon}(x)\exp\left[-\frac{2}{\epsilon}\left(\int_0^xF'\ast u_{\epsilon}(y)dy+V(x)\right)\right],
$$
with
$$\Lambda_{\epsilon}'(x)=\frac{2}{\epsilon}\, C_{\epsilon}\exp\left[\frac{2}{\epsilon}\left(\int_0^xF'\ast u_{\epsilon}(y)dy+V(x)\right)\right].$$
Hence
\begin{eqnarray*}
u_{\epsilon}(x)&=&\Lambda_{\epsilon}(0)\exp\left[-\frac{2}{\epsilon}\left(\int_0^xF'\ast u_{\epsilon}(y)dy+V(x)\right)\right]\\
&+&\frac{2}{\epsilon}\,C_{\epsilon}\int_0^x\exp\left[\frac{2}{\epsilon}\left(\int_0^yF'\ast u_{\epsilon}(z)dz+V(y)\right)\right]dy\\
&\times&\exp\left[-\frac{2}{\epsilon}\left(\int_0^xF'\ast u_{\epsilon}(y)dy+V(x)\right)\right].
\end{eqnarray*}
Let us assume that $C_{\epsilon}\neq0$. Applying Lemma \ref{lem:Annexe1} to the function $U(x)=\int_0^xF'\ast u_{\epsilon}(y)dy+V(x)$, whose second derivative is positive for $\vert x\vert$ large enough (using hypotheses (V-5) and (F-2)), permits to exhibit the equivalent of $\Lambda_{\epsilon}(x)$:
$$
\Lambda_{\epsilon}(x)\approx\frac{2}{\epsilon}C_{\epsilon}\frac{\exp\left[\frac{2}{\epsilon}\left(\int_0^xF'\ast u_{\epsilon}(y)dy+V(x)\right)\right]}{\frac{2}{\epsilon}\Big(V'(x)+F'\ast u_{\epsilon}(x)\Big)}\quad \mbox{as}\ x\to\pm\infty.
$$
Hence
$$u_{\epsilon}(x)\approx\frac{C_{\epsilon}}{V'(x)+F'\ast u_{\epsilon}(x)}.$$
Due to the conditions (V-6) and (F-4), there exists some constant $K>0$ such that
\[
\vert V'(x)+F'\ast u_{\epsilon}(x)\vert \le K(1+\vert x\vert ^{4q-1}),\quad \mbox{for all}\ x\in\mathbb{Re}.
\]
We deduce that $x\rightarrow x^{8q^2}$ can't be integrated with
respect to $u_{\epsilon}$: that contradicts the essential
assumption of the statement. We deduce that $C_{\epsilon}=0$ and
obtain \eqref{eq:lem:structure} after normalization.
\end{proof}
Lemma \ref{lem:exp_form} presents the essential structure of any invariant measure. The global exponential form will play a crucial role in next sections: to prove the existence of some stationary measure, it is necessary and sufficient to solve equation \eqref{eq:lem:structure}.
\section{The linear interaction case}
First we shall analyze the existence problem for stationary measures in the simple linear case. In this case $F'(x)=\alpha x$ with $\alpha>0$, the interaction gradient function is quadratic: $F(x)=\frac{\alpha}{2}\, x^2$ and the stochastic differential equation takes an interesting simple form. The non-linearity of the drift term is limited to the average of the density $u_t^\epsilon(x)$:
\[
X_t^{\epsilon}=X_0+\sqrt{\epsilon}B_t-\int_0^tV'(X_s^{\epsilon})ds-\alpha\int_0^t\Big( X_s^{\epsilon}-\int_{\mathbb{Re}}xdu_s^{\epsilon}(x)\Big)ds, \quad \epsilon>0 .
\]
The study of this particular case emphases the existence of several invariant measures. The interesting problem is then to determine in which situations the number of such measures is perfectly known.
\subsection{Existence of invariant measures}
The existence question is really simplified in the linear interaction case, it is just reduced {\itshape in fine} to the following parametrization problem. Let us denote the first moment of an invariant measure $u_{\epsilon}$ by
\begin{equation}
\label{eq:moyenne}
m_1(\epsilon)=\int_{\mathbb{Re}}xu_{\epsilon}(x)dx,
\end{equation}
then \eqref{eq:lem:structure} becomes
\begin{eqnarray}
\label{eq:structure:lineaire}
u_{\epsilon}(x)=\frac{\exp\left[-\frac{2}{\epsilon}\left(V(x)+\alpha\frac{x^2}{2}-\alpha m_1(\epsilon)x\right)\right]}{\int_{\mathbb{Re}}\exp\left[-\frac{2}{\epsilon}\left(V(y)+\alpha\frac{y^2}{2}-\alpha m_1(\epsilon)y\right)\right]dy}.
\end{eqnarray}
We now come to the essential equivalence: $u_{\epsilon}$ is an
invariant measure if and only if \eqref{eq:moyenne} and
\eqref{eq:structure:lineaire} are satisfied. It suffices then to
point out the convenient parameters $m_1(\epsilon)$ since there is
a one to one correspondence between these parameters and the
invariant measures. In other words, we shall find the solution of the equation 
\begin{equation}
\label{prop:proof:function}
m=\Psi_\epsilon(m)\quad\mbox{with}\ \Psi_{\epsilon}(m)=\frac{\int_{\mathbb{Re}}x\exp\left[-\frac{2}{\epsilon}\left(V(x)+\alpha\frac{x^2}{2}-\alpha mx\right)\right]dx}{\int_{\mathbb{Re}}\exp\left[-\frac{2}{\epsilon}\left(V(x)+\alpha\frac{x^2}{2}-\alpha mx\right)\right]dx}.
\end{equation}
Obviously, $m_1^0(\epsilon)=0$ is a
candidate. The corresponding measure $u_{\epsilon}^0$ is invariant
and symmetric:
\begin{eqnarray*}
 u_{\epsilon}^0(x)=\exp\left[-\frac{2}{\epsilon}\left(V(x)+\alpha\frac{x^2}{2}\right)\right]\left(\int_{\mathbb{Re}}\exp\left[-\frac{2}{\epsilon}\left(V(y)+\alpha\frac{y^2}{2}\right)\right]dy\right)^{-1}.
 \end{eqnarray*}
In fact $u_{\epsilon}^0$ is the unique \emph{symmetric} stationary measure.\\
Of course the natural question concerns the existence of others reals $m_1(\epsilon)$ solutions of \eqref{eq:structure:lineaire}. In fact the basic dynamical system associated to self-stabilizing diffusions is symmetric since $F$ and $V$ are assumed to be even functions. The consequence is immediate: if the initial law of the diffusion $(X_t^{\epsilon},\ t\ge 0)$ is symmetric so will be the law of $X_t^{\epsilon}$ for all $t>0$. In \cite{BRTV}, the authors consider self-stabilizing diffusions without the environment potential $V$. They proved the existence of some unique symmetric invariant measure and describe the behavior of the diffusion: for any initial law satisfying the moment condition of order $8q^2$ the law of $X_t-\mathbb{E}[X_0]$ converges to the invariant symmetric law as time elapses.\\
Adding some double-well potential $V$ in the main structure of the stochastic differential equation changes drastically the situation. In particular we prove the existence of several invariant measures, one of them being symmetric.
\begin{prop}
\label{prop:exist} Let $a$ be the unique positive real which
minimizes $V$ (see (V-3)). For all $\delta\in]0,1[$,
there exists $\epsilon_0>0$ such that for all
$\epsilon\leq\epsilon_0$, the equation
\eqref{prop:proof:function} admits a
solution satisfying the estimates:
\begin{eqnarray}
\label{eq:prop:exist}
a-\frac{(1+\delta)V^{(3)}(a)}{4V''(a)\left(\alpha+V''(a)\right)}\ \epsilon\leq m_1(\epsilon)\leq a-\frac{(1-\delta)V^{(3)}(a)}{4V''(a)\left(\alpha+V''(a)\right)}\ \epsilon.
\end{eqnarray}
Moreover $-m_1(\epsilon)$ satisfies \eqref{prop:proof:function} too.
\end{prop}
Let us note that, for $\epsilon$ small enough, the preceding
proposition implies the existence of at least three invariant
measures corresponding to the averages: $0$, $m_1(\epsilon)$ and
$-m_1(\epsilon)$.
\begin{proof} Set $\tau>0$. Let's proceed to the first
order asymptotic development of the expression
$\Psi_{\epsilon}(a-\tau\epsilon)$.
\begin{eqnarray*}
\Psi_{\epsilon}(a-\tau\epsilon)&=&\frac{\int_{\mathbb{Re}}x\exp\left[-\frac{2}{\epsilon}\left(V(x)+\alpha\frac{x^2}{2}-\alpha (a-\tau\epsilon)x\right)\right]dx}{\int_{\mathbb{Re}}\exp\left[-\frac{2}{\epsilon}\left(V(x)+\alpha\frac{x^2}{2}-\alpha (a-\tau\epsilon)x\right)\right]dx}\\
&=&\frac{\int_{\mathbb{Re}}xe^{-2\alpha\tau x}\exp\left[-\frac{2}{\epsilon}\left(V(x)+\alpha\frac{x^2}{2}-\alpha ax\right)\right]dx}{\int_{\mathbb{Re}}e^{-2\alpha\tau x}\exp\left[-\frac{2}{\epsilon}\left(V(x)+\alpha\frac{x^2}{2}-\alpha ax\right)\right]dx}.
\end{eqnarray*}
By Lemma \ref{lem:Annexe2} applied to the context:
$f(x)=-2\alpha\tau x$, $n=1$,
$U(x)=V(x)+\frac{\alpha}{2}x^2-\alpha x$ and $\mu=0$, we get:
\begin{eqnarray*}
\Psi_{\epsilon}(a-\tau\epsilon)&=&a-\frac{1}{4a\left(\alpha+V''(a)\right)^2}\left[aV^{(3)}(a)+4a\alpha\tau\left(\alpha+V''(a)\right)\right]\epsilon+o(\epsilon)\\
&=&a-\tau\epsilon+\frac{V''(a)}{\alpha+V''(a)}\left[\tau-\frac{V^{(3)}(a)}{4V''(a)\left(\alpha+V''(a)\right)}\right]\epsilon+o(\epsilon).
\end{eqnarray*}
Set $\tau^0=\frac{V^{(3)}(a)}{4V''(a)\left(\alpha+V''(a)\right)}$.
Then $a-\tau^0\epsilon$ is the first order approximation of the
fixed point. Indeed for $\delta\in]0;1[$ we
can define
\begin{eqnarray*}
d_{\pm}:=\Psi_{\epsilon}\left(a-\tau^0(1\pm\delta)\epsilon\right)-\left(a-\tau^0(1\pm\delta)\epsilon\right)
=\pm\delta\frac{V''(a)}{\alpha+V''(a)}\tau^0\epsilon+o(\epsilon).
\end{eqnarray*}
For $\epsilon$ small enough, $d_+>0$ and $d_-<0$. Since the
function $\Psi_{\epsilon}$ is $\mathcal{C}^0$ continuous, there
exists
$m_1(\epsilon)\in[a-\tau^0(1+\delta)\epsilon;\;a-\tau^0(1-\delta)\epsilon]$
which satisfies $\Psi_{\epsilon}(m_1(\epsilon))=m_1(\epsilon)$.
Finally, by the change of variable $x:=-x$ in the integral
expression \eqref{prop:proof:function}, we obtain
$\Psi_{\epsilon}(-m_1(\epsilon))=-\Psi_{\epsilon}(m_1(\epsilon))=-m_1(\epsilon)$.
\end{proof}
\subsection{Description of the set of invariant measures}\label{sect:set}
According to Proposition \ref{prop:exist}, we know there are at
least three invariants measures. One of them is symmetric
corresponding to the average $0$ and two others will be called
\emph{outlying measures}, one wrapped around $a$ and the other one
around $-a$. The aim of this section is to study if
there are exactly three invariants measures or more.\\
For this purpose, we study the asymptotic behavior of the function
$\Psi_{\epsilon}$ defined by \eqref{prop:proof:function} in the
small noise limit.
\begin{thm}
\label{prop:convex} If $V''$ is a convex function then, in the
small noise limit,  there exist exactly three stationary
measures.
\end{thm}
\begin{proof}
Let $m>0$. Let us recall that the interaction function is linear:
$F'(x)=\alpha x$ with $\alpha>0$. In order to study the invariant measures, we have to consider the fixed points of the application $\Psi_\epsilon(m)$ defined by \eqref{prop:proof:function}. We introduce the following potential function:
\[
W_m(x)=V(x)+\frac{\alpha}{2}x^2-\alpha mx.
\]
Since $V'(0)=0$, we have $W_m'(0)<0$. Moreover $\lim_{x\to+\infty}W_m'(x)=+\infty$. So we denote by $x_m$ the positive real for which the potential $W_m$ admits its global minimum. It is uniquely determined since $V''$ is a convex function. In particular, $x_m$ satisfies
 $V'(x_m)+\alpha(x_m-m)=0$ and $V''(x_m)+\alpha\geq0$. Furthermore $V''(x_m)+\alpha>0$. Indeed, since $x_m$ is a global minimum, the equality $V''(x_m)+\alpha=0$ implies that $V^{(3)}(x_m)=0$ that is $x_m=0$ which contradicts the assumption concerning the positivity of $x_m$. \\
We define
\[
\chi_{\epsilon}(m)=\Psi_{\epsilon}(m)-m\quad \mbox{and}\ \chi_0(m)=x_m-m.
\] 
We obtain the expression:
\begin{equation}\label{eq:proof:chi}
\chi_{\epsilon}(m)=x_m-m+\frac{\int_{\mathbb{R}}(x-x_m)\exp\left[-\frac{2}{\epsilon}\left(V(x)+\frac{\alpha x^2}{2}-\alpha mx\right)\right]dx}{\int_{\mathbb{R}}\exp\left[-\frac{2}{\epsilon}\left(V(x)+\frac{\alpha x^2}{2}-\alpha mx\right)\right]dx}.
\end{equation}
It suffices to prove that $\chi_\epsilon$ has just one zero in $\mathbb{Re}_+^*$. \\
{\bf Step 1: } {\it For all $\epsilon>0$ and $m>0$, we observe that $\chi_{\epsilon}(m)\le\chi_0(m)=x_m-m$.}\\
We apply the change of variable $x:=y+x_m$ to the integrals in \eqref{eq:proof:chi} and obtain
\begin{eqnarray*}
\chi_{\epsilon}(m)&=&\chi_0(m)+\frac{\int_{\mathbb{R}}y\exp\left[-\frac{2}{\epsilon}\left(V(y+x_m)+\frac{\alpha y^2}{2}+\alpha\left(x_m-m\right)y\right)\right]dy}{\int_{\mathbb{R}}\exp\left[-\frac{2}{\epsilon}\left(V(y+x_m)+\frac{\alpha y^2}{2}+\alpha\left(x_m-m\right)y\right)\right]dy}\\
&=&\chi_0(m)+\frac{\int_0^{\infty}y\exp\left[-\frac{\alpha}{\epsilon}y^2\right]\Omega_{\epsilon,m}(y)dy}{\int_{\mathbb{R}}\exp\left[-\frac{2}{\epsilon}\left(V(y+x_m)+\frac{\alpha y^2}{2}+\alpha\left(x_m-m\right)y\right)\right]dy},
\end{eqnarray*}
with
\begin{eqnarray*}
\Omega_{\epsilon,m}(y)&=&\exp\left[-\frac{2}{\epsilon}\Big(V(y+x_m)+\alpha\left(x_m-m\right)y\Big)\right]\\
&&-\exp\left[-\frac{2}{\epsilon}\Big(V(y-x_m)-\alpha\left(x_m-m\right)y\Big)\right].
\end{eqnarray*}
We introduce the function
\begin{eqnarray*}
\Lambda_m(y)=V(y+x_m)-V(y-x_m)+2\alpha\left(x_m-m\right)y
\end{eqnarray*}
Since $V$ is an even function, $\Lambda_m(0)=0$ and $\Lambda_m''(0)=0$. According to the definition of $x_m$, $\Lambda_m'(0)=0$. $V''$ is a convex function therefore $V^{(3)}$ is increasing. So $\Lambda_m^{(3)}(y)=V^{(3)}(y+x_m)-V^{(3)}(y-x_m)\ge 0$ for all $y$. We deduce that $\Lambda_m''$ is increasing. Hence $\Lambda_m''$ is nonnegative on $\mathbb{R}_+^{*}$ so does $\Lambda_m(y)$ for $y>0$. Finally we get $\Omega_{\epsilon,m}(y)
\le0$ for all $y>0$. We obtain the announced result: $\chi_{\epsilon}(m)\le\chi_0(m)$ for $m>0$.\\
{\bf Step 2.} {\it $\chi_0$ has a unique zero on $\mathbb{Re}_+^*$.}\\
Let us compute $\chi_0(a)$ with $a$ defined in (V-3). We know that $a$ is solution of  $V'(x)+\alpha\left(x-a\right)=0$ with $V''(x)+\alpha>0$. Hence $\chi_0(a)=0$.\\
Let us focus our attention to the variations of the function $\chi_0$ on the interval $]0,a]$. 
Since
$V'(x_m)+\alpha x_m=\alpha m,$
and $\alpha+V''(x_m)>0$ we deduce that $m\to x_m$ is derivable; we obtain 
$$\chi_0'(m)=\frac{d}{dm}\,x_m-1.$$
and
\begin{equation}\label{eq:proof:derivee}
\frac{d}{dm}\,x_m=\frac{\alpha}{\alpha+V''(x_m)}>0\quad\mbox{which implies}\ \chi_0'(m)=-\frac{V''(x_m)}{\alpha+V''(x_m)}.
\end{equation}
The denominator is positive due to the definition of $x_m$.
According to (V~-~5), $V''(x)>0$ for all $x>a$. Hence $\chi_0'(m)<0$ for all $m>a$. Since $\chi_0(a)=0$ we deduce that, for all $m>a$,  $\chi_0(m)$ is strictly negative and therefore the function $\chi_{0}$ has no zero on $]a;+\infty[$.\\
It remains to study $\chi_0$ on the interval $]0,a]$. Since $V''$ is a convex function, we deduce that the derivative of $\chi_0$ is non positive for $x_m\ge c$ with $c>0$ satisfying $V''(c)=0$. We know that $c>0$ is unique since $V''(0)<0$ and $V''$ is a convex function. Moreover $c<a$. Since the function $m\to x_m$ is increasing for $m>0$, we deduce that $\chi_0'$ is negative for $x\in]\max(0,m_c),a]$ where $m_c=c+\frac{V'(c)}{\alpha}$. By construction, if $m_c>0$ then the equality $x_{m_c}=c$ holds.\\
We observe then two different cases:
\begin{itemize}
\item If $m_c\le 0$ i.e. $\alpha<\frac{\left|V'(c)\right|}{c}$: $\chi_0$ is decreasing on $\mathbb{Re}_+^*$ with $\chi_0(a)=0$. The unique zero of $\chi_0$ on $\mathbb{Re}_+^*$ is $a$.
\item If $m_c>0$ then $\chi_0$, which is a continuous function on $\mathbb{Re}_+^*$, is increasing on $]0,m_c[$ and decreasing on $]m_c,+\infty[$ with $\chi_0(a)=0$. It suffices to prove that $\lim_{m\to 0+}\chi_0(m)\ge 0$ in order to conclude that $a$ is the unique zero of $\chi_0$ on $\mathbb{Re}_+^*$. Due to the definition of $x_m$ we get: $\lim_{m\to 0+}\chi_0(m)=\lim_{m\to 0+}x_m\ge 0$. Indeed $m\to x_m$ is continuous from $]0,+\infty[$ to $]0,+\infty[$ so the extension to $m=0$ is non negative. 
\end{itemize}
In these two cases, there is a unique zero of $\chi_0$ on $\mathbb{Re}_+^*$.\\
{\bf Step 3.} {\it The family of functions $(\chi_\epsilon)_\epsilon$ (respectively $(\chi_\epsilon')_\epsilon$) converges uniformly towards $\chi_0$ (resp. $\chi_0'$) on each compact subset of $\mathbb{Re}_+^*$.}\\ 
First we prove the convergence of $\chi_\epsilon(m)$ for $m>0$. Recall that 
\begin{eqnarray*}
\chi_{\epsilon}(m)=\frac{\int_{\mathbb{R}}x\exp\left[-\frac{2}{\epsilon}\left(V(x)+\frac{\alpha x^2}{2}-\alpha m x\right)\right]dx}{\int_{\mathbb{R}}\exp\left[-\frac{2}{\epsilon}\left(V(x)+\frac{\alpha x^2}{2}-\alpha mx\right)\right]dx}.
\end{eqnarray*}
By Lemma \ref{lem:Annexe2} with $U(x)=V(x)+\frac{\alpha x^2}{2}$, $n=1$, $\mu=m$ and $G=-\alpha x$ we obtain the announced convergence result:
\begin{eqnarray*}
\chi_{\epsilon}(m)-\chi_0(m)=\chi_\epsilon(m)-x_m+m=-\frac{V^{(3)}(x_m)}{4\left(\alpha+V''(x_m)\right)^2}\epsilon+o(\epsilon).
\end{eqnarray*}
Moreover this convergence is uniform with respect to the variable $m$ on compact subsets of $\mathbb{Re}_+^*$.\\
We estimate now the asymptotics of $\chi_\epsilon'(m)$ as $\epsilon$ becomes small. Taking the derivative of $\Psi_\epsilon$, we obtain
\begin{eqnarray*}
\Psi_{\epsilon}'(m)&=&\frac{2\alpha}{\epsilon}\left\{\frac{\int_{\mathbb{Re}}x^2\exp\left[-\frac{2}{\epsilon}W_m(x)\right]dx}{\int_{\mathbb{Re}}\exp\left[-\frac{2}{\epsilon}W_m(x)\right]dx}-\left(\frac{\int_{\mathbb{Re}}x\exp\left[-\frac{2}{\epsilon}W_m(x)\right]dx}{\int_{\mathbb{Re}}\exp\left[-\frac{2}{\epsilon}W_m(x)\right]dx}\right)^2\right\}.
\end{eqnarray*}
We recognize the variance of the measure $u_{\epsilon}^{(m)}$ which is the measure associated to the average $m$ by \eqref{eq:structure:lineaire}. Hence
\begin{equation}\label{eq:var2}
\chi_{\epsilon}'(m)=\frac{2\alpha}{\epsilon}{\rm Var}(u_{\epsilon}^{(m)})-1.
\end{equation}
Applying again Lemma \ref{lem:Annexe2} with $U=V(x)+\frac{\alpha x^2}{2}$, $G=-\alpha x$, $\mu=m$ and $n=2$, we obtain
\begin{eqnarray*}
\frac{\int_{\mathbb{Re}}x^2\exp\left[-\frac{2}{\epsilon}W_m(x)\right]dx}{\int_{\mathbb{Re}}\exp\left[-\frac{2}{\epsilon}W_m(x)\right]dx}=x_m^2-\frac{\left(x_mV^{(3)}(x_m)-(\alpha+V''(x_m))\right)}{2\left(\alpha+V''(x_m)\right)^2}\,\epsilon+o(\epsilon).
\end{eqnarray*}
Applying the same lemma with $n=1$ permits to compute the first moment:
\begin{eqnarray*}
\frac{\int_{\mathbb{Re}}x\exp\left[-\frac{2}{\epsilon}W_m(x)\right]dx}{\int_{\mathbb{Re}}\exp\left[-\frac{2}{\epsilon}W_m(x)\right]dx}=x_m-\frac{x_mV^{(3)}(x_m)}{4x_m\left(\alpha+V''(x_m)\right)^2}\epsilon+o(\epsilon).
\end{eqnarray*}
By \eqref{eq:var2} and the computations of the two first moments, we get
\begin{equation}\label{eq:result-import}
\chi_{\epsilon}'(m)=\frac{-V''(x_m)}{\alpha+V''(x_m)}+o(1)=\chi_0'(m)+o(1).
\end{equation}
Furthermore this convergence is uniform with respect to the variable $m$ on compact subsets of $\mathbb{Re}_+^*$.\\
{\bf Step 4.} {\it For any $\delta>0$ small enough, there exists $\epsilon_0>0$ such that $\chi_\epsilon$ has a unique zero on $[\delta,\infty[$ for all $\epsilon\le \epsilon_0$.}\\
Since there is no zero of $\chi_\epsilon$ on the interval $]a,+\infty[$ (Step 1 and 2), we focus our attention to the interval $]0,a]$. On each compact subset of this interval, $\chi_\epsilon$ converges uniformly towards the limit function $\chi_0$ (Step 3). Hence the zeros of $\chi_\epsilon$ are in a small neighborhood of the unique zero of $\chi_0$ namely $a$ (Step 2). Let us study the derivative of $\chi_\epsilon$ in a neighborhood of $a$.  Since $\chi'_\epsilon$ converges uniformly towards $\chi'_0$ (Step 3) and $\chi_0'(m)<0$ in a neighborhood of $a$ (Step 2), we obtain that $\chi_\epsilon'(m)<0$ in a neighborhood of $a$ for $\epsilon$ small enough. 
Finally we proved that, as soon as $\epsilon$ is small enough, the function $\chi_{\epsilon}$ can't admit two zeros or more on $\mathbb{R}_+^{*}$. \\
{\bf Step 5.} {\it There exists $\delta>0$ and $\epsilon_0>0$ such that $\chi_\epsilon$ doesn't vanish on $]0,\delta]$ for all $\epsilon\le\epsilon_0$.}\\
In this last step, we have to distinguish three different cases depending on the values $\theta$ and $\alpha$ defined by \eqref{eq:param_theta} and \eqref{def:alpha}.\\
{\bf Step 5.1.} We assume $\alpha<\theta$. In this particular case $W_0(x)=V(x)+\alpha x^2/2$ reaches a unique global minimum on $\mathbb{Re}_+$ for $x=x_0>0$.\\
Let us fix some small $\delta>0$ (depending on $x_0$: we shall precise it in the following). We prove that, for $\epsilon$ small enough, $\chi_\epsilon(m)=\Psi_\epsilon(m)-m>0$ on $]0,\delta]$. By the definition of $\Psi_\epsilon$, see \eqref{prop:proof:function}, it suffices to prove that $N_\epsilon(m)>0$ for $m\in]0,\delta]$ where 
\begin{equation}
\label{def:n_eps}
N_\epsilon(m)=\int_\mathbb{Re} x\exp\Big[ -\frac{2}{\epsilon}\,W_m(x) \Big]dx-m\int_\mathbb{Re}\exp\Big[ -\frac{2}{\epsilon}\,W_m(x) \Big]dx.
\end{equation}
Obviously $N_\epsilon(0)=0$. Let us prove that $N_\epsilon$ is non decreasing.
Taking the derivative, we get
\[
N_\epsilon'(m)=\frac{2\alpha}{\epsilon}\int_\mathbb{Re}\Big( x^2-mx-\frac{\epsilon}{2\alpha} \Big)\exp\Big[ -\frac{2}{\epsilon}\,W_m(x) \Big]dx.
\]
This expression is in fact non negative. Indeed, using the symmetry property of $W_0(x)$ and the upper bound $m\le\delta$, we obtain
\begin{align*}
N_\epsilon'(m)&=\frac{2\alpha}{\epsilon}\int_0^\infty\left\{ \Big( x^2-\frac{\epsilon}{2\alpha} \Big)\cosh\Big( \frac{2\alpha mx}{\epsilon} \Big)-mx\sinh\Big( \frac{2\alpha mx}{\epsilon} \Big) \right\}e^{-\frac{2}{\epsilon}\,W_0(x)}dx\\
&\ge \frac{\alpha}{\epsilon}\int_0^\infty P_\delta(x) e^{\frac{2\alpha mx}{\epsilon}}e^{-\frac{2}{\epsilon}\,W_0(x)}dx\quad\mbox{with}\ P_\delta(x)=x^2-\delta x-\frac{\epsilon}{\alpha}.
\end{align*}
We split the preceding integral into two parts: the first integral $I_0$ concerns the support $[0,2\delta]$ and the second integral $I_{2\delta}$ the complementary support $[2\delta,\infty[$. We get $N_\epsilon'(m)\ge \frac{\alpha}{\epsilon}\, (I_0+I_{2\delta})$.\\
Since the roots of the polynomial function $P_\delta$ satisfy 
\[
x_\pm=\frac{1}{2}\Big( \delta\pm\sqrt{\delta^2+\frac{4\epsilon}{\alpha}} \Big)<2\delta,
\]
the polynomial is positive on the interval $[2\delta,\infty[$ and can be lower bounded by $P_\delta(2\delta)=2\delta^2-\epsilon/\alpha$. Lemma \ref{lem:Annexe:inter2} implies the existence of some constant $C>0$ leading to the following estimate as $\epsilon\to 0$:
\begin{align}\label{premm}
I_{2\delta}\ge (2\delta^2-\epsilon/\alpha) \int_{2\delta}^{x_0+1}e^{-\frac{2}{\epsilon}\,W_0(x)}dx\ge C\delta^2\sqrt{\epsilon}\ e^{-\frac{2}{\epsilon}(V(x_0)+\alpha x_0^2/2)}
\end{align}
provided that $x_0>2\delta$ (it suffices then to chose $\delta$ small enough).\\
Let us finally focus our attention to the lower bound of the integral term $I_0$. Since the minimum value of $P_\delta$ is $-(\delta^2/4+\epsilon/\alpha)$ and since $W''(0)<0$, we have
\begin{align}\label{premmm}
I_0\ge -\Big(\frac{\delta^2}{4}+\frac{\epsilon}{\alpha} \Big)\int_0^{2\delta}e^{\frac{2\alpha mx}{\epsilon}}e^{-\frac{2}{\epsilon}\,W_0(x)}dx\ge -2\delta\Big(\frac{\delta^2}{4}+\frac{\epsilon}{\alpha} \Big) e^{-\frac{V(2\delta)}{\epsilon}}.
\end{align}
For $\delta>0$ small enough, $V(2\delta)>V(x_0)+\alpha x_0^2/2$ (since the minimum of $V(x)+\alpha x^2/2$ is only reached for $x=x_0$). Consequently the negative lower bound of $I_0$ \eqref{premmm} is negligible with respect to the positive lower bound of $I_{2\delta}$ as $\epsilon$ becomes small. We deduce that there exists $\epsilon_0$ such that $N_\epsilon'(m)>0$ for all $m\in[0,\delta]$ and $\epsilon\le\epsilon_0$. Since $N_\epsilon(0)=0$ we conclude that $N_\epsilon(m)>0$ on $]0,\delta]$ and so is $\chi_\epsilon$.\\
{\bf Step 5.2.} We assume $\alpha>\theta$. In this case $W_0(x)$ admits a unique minimum reached for $x=0$ and $x_m$ converges continuously to $0$ as $m\to 0$. Using similar arguments as those presented in Step 3, we claim that $\chi_\epsilon$ (resp. $\chi'_\epsilon$) converges towards $\chi_0$ (resp. $\chi_0'$) uniformly on $[0,a]$ as $\epsilon\to 0$. 
Due to the regularity of $\chi_0$ and by the inequality $\chi_0'(0)=-\frac{V''(0)}{\alpha+V''(0)}>0$ we obtain the existence of $\delta>0$ and $\epsilon_0>0$ such that $\chi'_\epsilon(m)>0$ for $m\in[0,\delta]$ and $\epsilon\le\epsilon_0$. $\chi_\epsilon$ starts in $0$ and is strictly increasing on $[0,\delta]$ which implies the announced result.\\
{\bf Step 5.3.} We assume that $\alpha=\theta$. It suffices then to note that $\chi_\epsilon$ depends continuously on the parameter $\alpha$. The following results can be directly deduced from the preceding case (Step 5.2) by continuity: $\chi_\epsilon'(0)>0$ and $\chi_\epsilon'(m)\ge 0$ for $m\in[0,\delta]$ and $\epsilon\le\epsilon_0$. In fact $\chi_\epsilon$ vanishes for  $x=0$ and is increasing on $[0,\delta]$. The inequality $\chi_\epsilon(m)>0$ for all $m\in]0,\delta]$ and $\epsilon\le\epsilon_0$ is an obvious consequence. \\
{\bf Conclusion:} Step 4 and 5 lead to the existence of $\epsilon_0>0$ such that for all
$\epsilon<\epsilon_0$, $\chi_{\epsilon}$ has exactly three zeros:
$0$ and two other reals, one in the neighborhood of $a$, the other one near $-a$. To each of these averages corresponds a unique invariant measure obtained by \eqref{eq:structure:lineaire}.
\end{proof}
\noindent{\bf Example:} In Theorem \ref{prop:convex}, for all $\alpha>0$, as soon as $\epsilon$ is small enough, there exist exactly three invariant measures. There is a one to one correspondence between these measures and their average through \eqref{eq:structure:lineaire}. It suffices to determine the averages which are in fact solutions to the equation 
\[
\chi_\epsilon^\alpha(m):=\Psi_\epsilon(m)-m=0.
\]
These solutions are really close to the solutions of $\chi_0^\alpha(m)=x_m^\alpha-m=0$ in the small noise limit. We recall that $x_m^\alpha$ is the global minimum of
\[
W_m^\alpha(x):=V(x)+\frac{\alpha}{2}x^2-\alpha mx\quad \mbox{on}\ \mathbb{Re}_+^*.
\] 
Let us observe these averages in the particular case:
$V(x)=\frac{x^4}{4}-\frac{x^2}{2}$ and $F(x)=\frac{\alpha}{2}x^2$. In this case, we compute the parameter $c>0$ which vanishes $V''$ and the corresponding parameter $m_c=c+\frac{V'(c)}{\alpha}$. We obtain:
\begin{equation}\label{param}
c=\frac{1}{\sqrt{3}}\quad\mbox{and}\ m_c=\frac{3\alpha-2}{3\sqrt{3}\alpha}.
\end{equation}
We shall for this example present graphs of the functions $\chi_0^\alpha$ (dotted line) and $\chi_\epsilon^\alpha$ for different values of $\alpha$. We choose $\epsilon=1/4$. Even if it seems to be not very small, this value suffices in this example to observe three invariant measures for each interaction parameter value considered.

First of all we have to determine the value of $x_m^\alpha$ which is solution of the system $\left(E^{\alpha,m}\right)$:
\begin{eqnarray*}
X^3+(\alpha-1)X-\alpha m=0\quad\mbox{and}\quad3X^2+(\alpha-1)\geq0.
\end{eqnarray*}
Its discriminant is equal to
$$
\Delta_{\alpha}(m)=\frac{\alpha^2m^2}{4}+\frac{(\alpha-1)^3}{27}.
$$
We distinguish different cases:
\begin{itemize}
\item $\alpha=0$: the solution is evident, we get $x_m^\alpha=1$ and $\chi_0^0(m)=1-m$ for $m>0$ and by symmetry $\chi_0^0(m)=-1-m$ for all $m<0$. Moreover $\chi_0^0(0)=0$.
\item $\alpha>1$ (Figure \ref{figure>1}): for all $m\in\mathbb{R}$, we get $\Delta_{\alpha}(m)>0$. Hence
    \begin{eqnarray*}
    \chi_0^{\alpha}(m)=\sqrt[3]{\frac{\alpha m}{2}+\sqrt{\Delta_\alpha(m)}}+\sqrt[3]{\frac{\alpha m}{2}-\sqrt{\Delta_\alpha(m)}}-m.
    \end{eqnarray*}
    The function $\chi_0^{\alpha}$ is $\mathcal{C}^\infty$-continuous and odd. We observe also that $\chi_0^{\alpha}(0)=\chi_0^{\alpha}(1)=0$ and ${\chi_0^\alpha}'(m_c)=0$ with $m_c$ defined by \eqref{param}. Hence $\chi_0^\alpha$ is increasing on $]0,m_c[$ and decreasing on $]m_c,\infty[$. 
\item $\alpha=1$ (Figure \ref{figure3}) then $\chi_0^\alpha(m)=m^{\frac{1}{3}}-m$. The limit function is odd, continuous on $\mathbb{R}$ and $\mathcal{C}^{\infty}$ on $\mathbb{Re}^*$. Moreover the path is increasing for $m\in]0,m_c[$, decreasing for $m\in]m_c,\infty[$ with $m_c=\frac{1}{3\sqrt{3}}$.
\item $\frac{2}{3}<\alpha<1$ (Figure \ref{figure4}): the discriminant can be negative. Therefore let us define $m_0(\alpha)$ such that $\Delta_{\alpha}(m_0(\alpha))=0$. Then for all $m$ between $0$ and $m_0(\alpha)$, the discriminant is negative and for all $m$ larger than $m_0(\alpha)$ it is positive. We get $m_0(\alpha)=2(1-\alpha)^{\frac{3}{2}}/(3\alpha\sqrt{3})$.
We obtain the following function: $\chi_0^\alpha(0)=0$ and
\begin{eqnarray*}
\chi_0^{\alpha}(m)=\left\{
\begin{array}{lcr}
\varphi_1^{(\alpha)}(m)\quad\forall m\in[-m_0(\alpha);0[\bigcup]0;m_0(\alpha)]\\
\varphi_2^{(\alpha)}(m)\quad\forall m\in]-\infty;-m_0(\alpha)]\bigcup[m_0(\alpha);+\infty[
\end{array}
\right.
\end{eqnarray*}
    with
    \begin{eqnarray*}
    \varphi_1^{(\alpha)}(m)&=&2\sqrt{\frac{1-\alpha}{3}}\cos\left[\frac{1}{3}\arccos\left(\frac{\alpha m}{2}\sqrt{\frac{27}{(1-\alpha)^3}}\right)\right]-m\\
    \varphi_2^{(\alpha)}(m)&=&\sqrt[3]{\frac{\alpha m}{2}+\sqrt{\Delta_\alpha(m)}}+\sqrt[3]{\frac{\alpha m}{2}-\sqrt{\Delta_\alpha(m)}}-m.
    \end{eqnarray*}
Let us note that $\chi_0^\alpha(0^+)=\sqrt{1-\alpha}\neq0$ and $\chi_0^\alpha(0^-)=-\sqrt{1-\alpha}\neq0$. The function is $\mathcal{C}^{\infty}$-continuous on $]0;m_0(\alpha)[\cup]m_0(\alpha);+\infty[$ and continuous in $m_0(\alpha)$.\\
Moreover the function is increasing on the interval $]0,m_c[$ and decreasing for $m>m_c$. The maximum is therefore reached for $m=m_c$. We observe that $m_c\le m_0^\alpha$ for $\alpha\in[2/3,3/4]$ and $m_c\ge m_0^\alpha$ for $\alpha\in[3/4,1]$. We remark also that the increasing part is smaller and the decreasing part is longer for smaller values of $\alpha$.\\
Furthermore, the part where $\chi_{\alpha}$ is equal to $\varphi_1^{(\alpha)}$ is longer.
\item $\alpha\le \frac{2}{3}$ (Figure \ref{figure5}): the function $\chi_0^\alpha$ is defined in the same way as in the preceding case. The important difference is that the function is decreasing on $\mathbb{Re}_+^*$ since $m_c$  defined by \eqref{param} is non positive. 
\end{itemize}
\begin{minipage}[l]{5.5cm}
\begin{figure}[H]
\centerline{\includegraphics[width=4cm,angle=270]{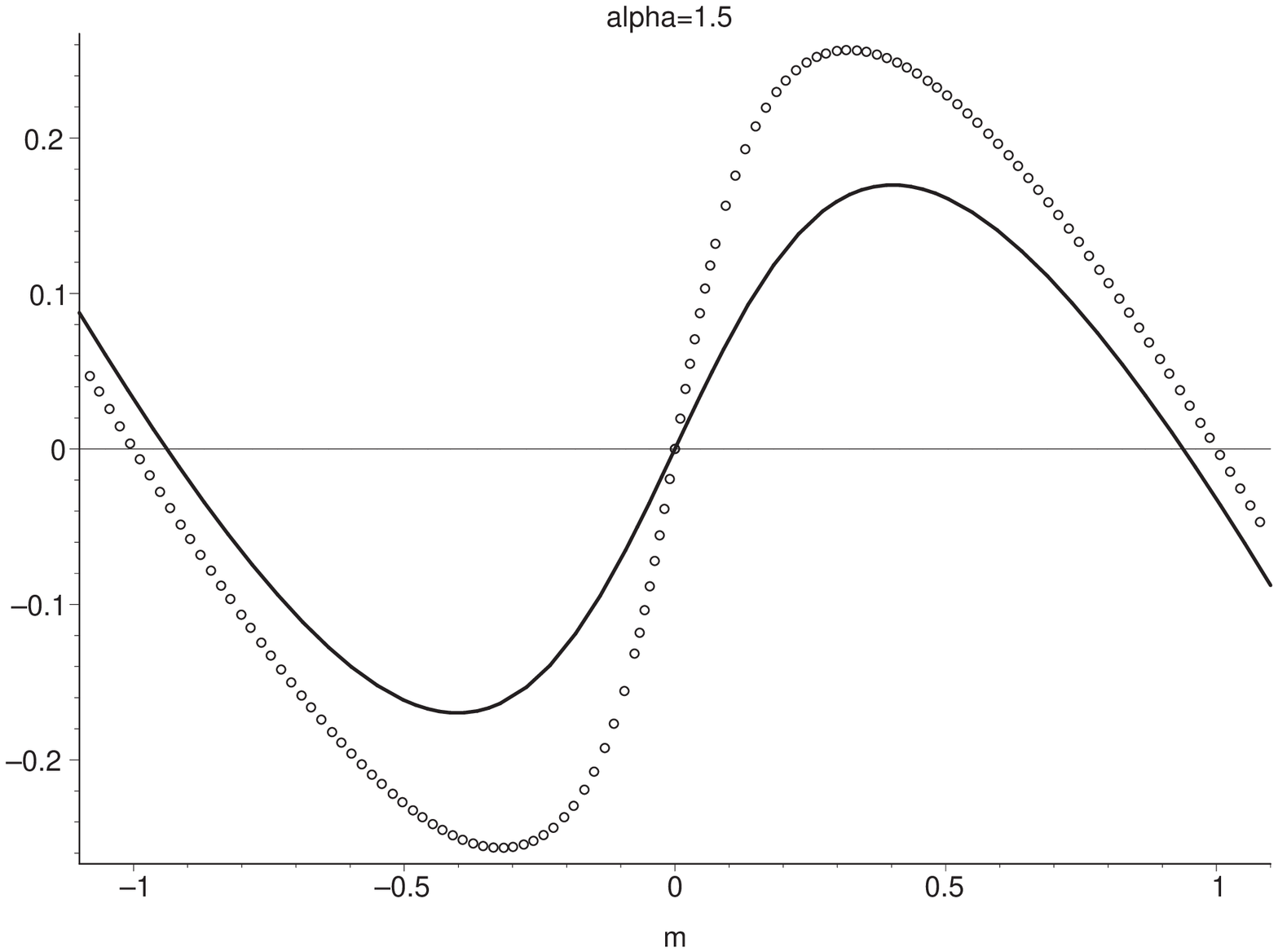}}
\caption{$\chi_0^\alpha$ (dotted line) and $\chi_\epsilon^\alpha$ for $\alpha>1$}\label{figure>1}
\end{figure}
\end{minipage}\ \ \ 
\begin{minipage}[l]{5.5cm}
\begin{figure}[H]
\centerline{\includegraphics[width=4cm,angle=270]{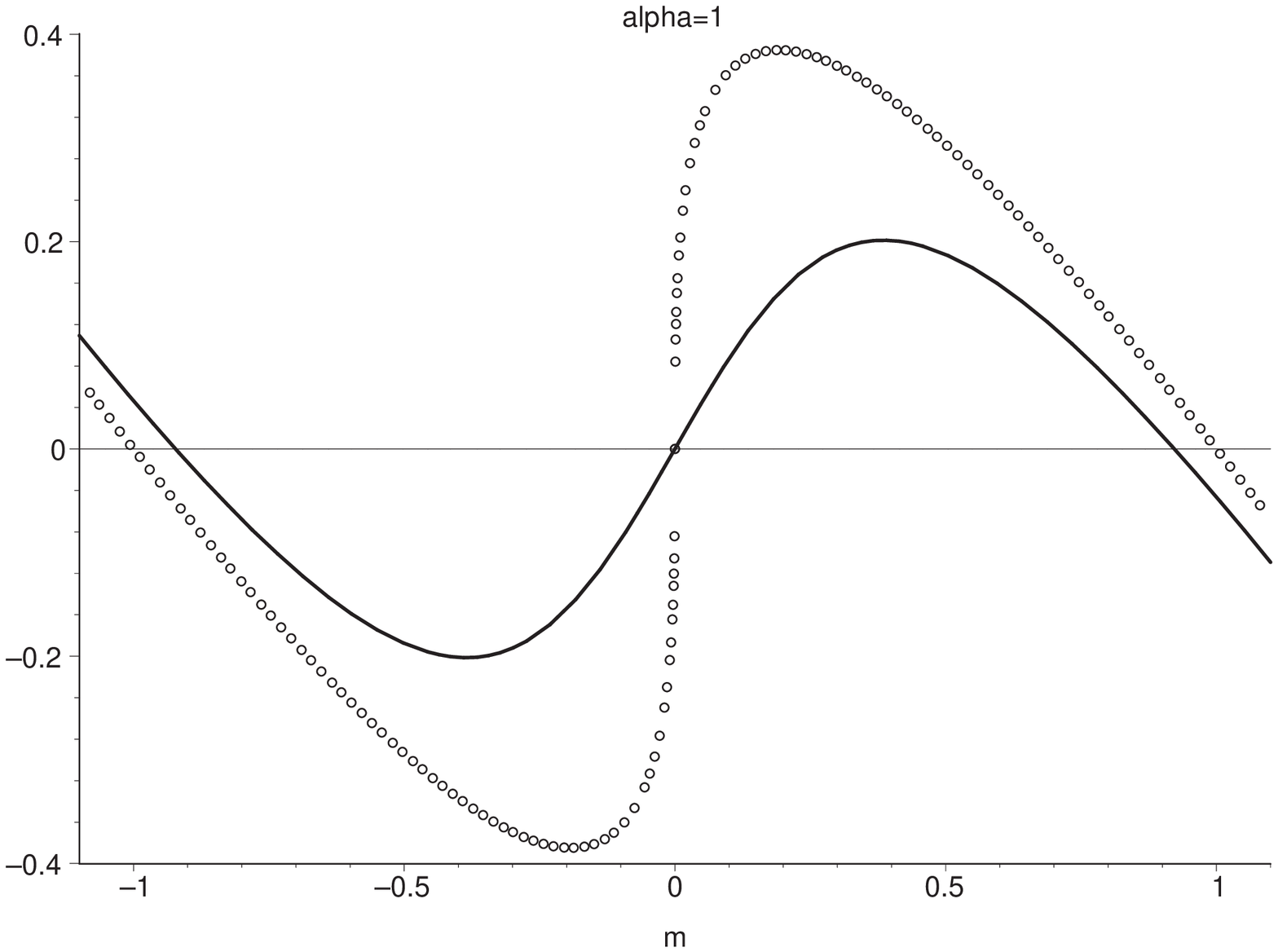}}
\caption{$\chi_0^\alpha$ (dotted line) and $\chi_\epsilon^\alpha$ for $\alpha=1$ }\label{figure3}
\end{figure}
\end{minipage}\\
\begin{minipage}[l]{5.5cm}
\begin{figure}[H]
\centerline{\includegraphics[width=4cm,angle=270]{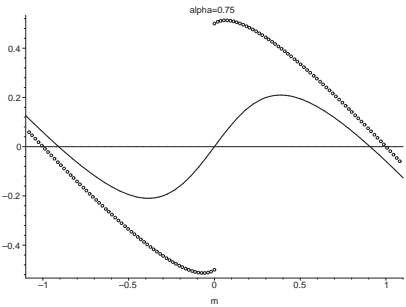}}
\caption{$2/3< \alpha<1$}\label{figure4}
\end{figure}
\end{minipage}\ \ \ 
\begin{minipage}[l]{5.5cm}
\begin{figure}[H]
\centerline{\includegraphics[width=4cm,angle=270]{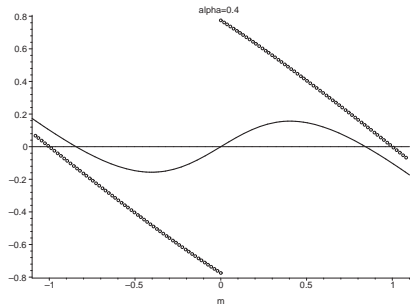}}
\caption{$\alpha<2/3$}\label{figure5}\end{figure}
\end{minipage}
\section{The general interaction case}
We assumed for this study that the self-attraction phenomenon is represented by a polynomial function $F'$, see (F-1). In previous section, we analyzed the particular linear situation: $F'(x)=\alpha x$ and proved under suitable conditions that there exist exactly three invariant measures in the small noise limit. In this section we shall focus our attention to the general case: the polynomial function $F$ is of degree $n\ge 2$. First we shall present results concerning the symmetric invariant measure and secondly we discuss the presence of asymmetric measures.
\subsection{Symmetric invariant measures}\label{sym}
In the linear case we proved the existence of a unique symmetric invariant measure. The result is obvious since it suffices  to solve the equation \eqref{prop:proof:function} with $m_1(\epsilon)=0$. In the general case in order to find the symmetric measure we have to solve some equation like \eqref{eq:structure:lineaire} but depending on much more parameters than just the mean $m_1(\epsilon)$. The total number of parameters depends in fact on the degree of $F$. Instead of trying to solve such system, we choose some other kind of proof based on a fixed point theorem which permits to prove the existence of symmetric invariant measures in even more general cases: the interaction function does not need to be polynomial. In \cite{BRTV}, Benachour, Roynette, Talay and Vallois introduced this method of proof for a self-stabilizing diffusion in the constant environment case ($V'(x)=0$). This proof can be adapted to our situation and is based on the following Schauder's theorem (see for instance \cite{Gilbarg77} Corollary 11.2 p. 280):
\begin{prop}\label{schauder}
Let $\mathbb{B}$ a Banach space, $\mathbb{C}$ a closed convex subset and $\mathbb{A}$ a continuous application $\mathbb{C}\to\mathbb{C}$ such that $\overline{\mathbb{A}(\mathbb{C})}$ is compact.
Then $\mathbb{A}$ admits a fixed point in $\mathbb{C}$.
\end{prop}
In order to use this proposition we introduce some definitions and notations:
\begin{enumerate}
  \item Let us choose $p>4q$ where $q$ is defined in (V-6).
  \item $\mathbb{D}=\left\{v:\mathbb{Re}\longrightarrow\mathbb{Re}^+\ \left|\right.\ v\ \mbox{is symmetric and}\  \sup_{x\in\mathbb{Re}^+}\left(1+|x|^{p}\right)v(x)<\infty\right\}$.
    \item $\mathbb{B}=\{f:\mathbb{Re}\longrightarrow\mathbb{Re}\ ;\ \sup_{x\in\mathbb{Re}}\left(1+|x|^p\right)|f(x)|<\infty\}$. Let us note that $\mathbb{D}\subset\mathbb{B}$. $\mathbb{B}$ is equipped with the norm $|\cdot |_{\infty}$ where $|f|_{\infty}=\sup_{x\in\mathbb{Re}}\left(1+|x|^p\right)|f(x)|$.
    \item For all $M>0$ we define the function space $\mathbb{C}_{M}$ as the subset of all non negative and even function belonging to $\mathbb{B}$ which satisfy:
\begin{eqnarray*}
\int_{\mathbb{Re}} f(x)\,dx=1\quad\mbox{and}\ \sup_{x\in\mathbb{Re}}\left(1+|x|^p\right)f(x)\leq M.
\end{eqnarray*}
    \item For any function $f\in\mathbb{D}$ we define the operator:
\begin{eqnarray}\label{def:operateur:A}
\mathbb{A}^{\epsilon}(f)(x)&=&\frac{\exp\left[-\frac{2}{\epsilon}\left(V(x)+\int_0^x\left(F'\ast f\right)(y)dy\right)\right]}{\int_{z\in\mathbb{Re}}\exp\left[-\frac{2}{\epsilon}\left(V(z)+\int_0^z\left(F'\ast f\right)(y)dy\right)\right]}\\
&=&\frac{1}{\lambda_{\epsilon}(f)}\exp\left[-\frac{2}{\epsilon}\left(V(x)+\int_0^x\left(F'\ast f\right)(y)dy\right)\right],\nonumber
\end{eqnarray}
where $\lambda_\epsilon(f)$ is the normalization factor.
    \item For any function $u\in\mathbb{D}$, we define the moments $\gamma_k(u)=\int_{\mathbb{Re}}|x|^ku(x)dx$ with $0\leq k\leq p-2$.
\end{enumerate}
Let us just point out that $\mathbb{C}_{M}$ is a closed and convex subset of $\mathbb{B}$. Moreover we have $\mathbb{C}_{M}\subset\mathbb{D}\subset\mathbb{B}$. The aim of this section will consist in proving that the application $\mathbb{A}^{\epsilon}$ is $\mathcal{C}^0(\mathbb{C}_M,\mathbb{C}_M)$-continuous and that $\overline{\mathbb{A}^\epsilon(\mathbb{C}_M)}$ is compact. Therefore Schauder's theorem implies the existence of a fixed point and as the matter of fact the existence of an invariant measure in the function space $\mathbb{D}$.
\begin{lem}\label{lem:prem:major}
For all $u\in\mathbb{C}_{M}$, we have:
\begin{enumerate}
    \item $\gamma_k(u)\leq MC_1$ where $\displaystyle C_1=1+\max_{0\leq r\leq p-2}\int_{\mathbb{Re}}\frac{|x|^r}{1+|x|^p}dx$ for all $0\leq k\leq p-2$.
    \item there exists a constant $C_2>0$ independent of $M$ such that 
    \begin{equation}
    \label{eq:lem:prem:major}
    \frac{\alpha}{2}x^2\leq\int_0^x(F'\ast u)(y)dy\leq C_2Mx^2(1+x^{2q})\quad\mbox{for all}\ x.
    \end{equation}
\end{enumerate}\end{lem}
\begin{proof}
{\bf 1.} Let $u\in\mathbb{C}_{M}$ then the function $x\to\frac{|x|^k}{1+|x|^{p}}$ is integrable on $\mathbb{Re}$  since $k\leq p-2$. Moreover the definition of $\mathbb{C}_M$ implies that $\left(1+|x|^p\right)u(x)\leq M$ for all $x\in\mathbb{Re}$. Therefore
  $$
  \gamma_k(u)=\int_{\mathbb{Re}}\frac{|x|^k}{1+|x|^p}(1+|x|^p)u(x)dx\leq M\int_{\mathbb{Re}}\frac{|x|^k}{1+|x|^p}dx\leq MC_1.
  $$
{\bf 2.} Let $x\ge 0$. Since $u\in\mathbb{C}_{M}$, $u$ is an even function. By \eqref{def:alpha} we have $F'(x)=\alpha x+F'_0(x)$ and (F-3) implies that $F'$ and $F_0'$ are non negative odd functions so is $ F'_0\ast u$. Using the inequality developed in the statement of Lemma 4.3 in \cite{BRTV} and the assumption (F-3), we have
  \[
  F_0'(x)\le \frac{1}{2}\,\Big( F'_0(x-y)+F'_0(x+y) \Big)\quad\mbox{for} \ y\in\mathbb{Re},\ x\ge 0.
  \] 
Therefore, for $x\geq 0$:
\begin{eqnarray*}
\int_0^x\left(F'_0\ast u\right)(y)dy&=&\int_0^x\int_0^{\infty}\left(F'_0(y-z)+F'_0(y+z)\right)u(z)dzdy\\
&\geq&\int_0^x\int_{0}^\infty 2F'_0(y)u(z)dzdy\geq0
\end{eqnarray*}
From the preceding inequality we deduce
\[
\int_0^x\left(F'\ast u\right)(y)dy=\int_0^x\left(F'_0\ast u\right)(y)dy+\frac{\alpha}{2}x^2\geq\frac{\alpha}{2}x^2\quad\mbox{for all}\ x\geq 0.
\]
Since $\int_0^x (F'\ast u)(y)dy$ is an even function, we get the inequality for all $x\in\mathbb{Re}$.\\
{\bf 3.} Due to the symmetry of $F'\ast u$ we restrict our study to $x\geq 0$. 
  \[
  \int_0^x\left(F'\ast u\right)(y) dy=\frac{1}{2}\int_0^x\int_{0}^\infty \Big(F'(y-z)+F'(y+z)\Big)u(z)dz dy.
  \]
  According to the assumptions (F-1) and (F-4), $F$ is an even polynomial function of degree smaller than $2q$ with $q\geq1$. We can therefore write $F'$ as follows
  $$
  F'(x)=\sum_{k=0}^{q-1}\alpha_k x^{2k+1}.
  $$
  Therefore defining $\mathcal{F}(y,z)=F'(y-z)+F'(y+z)$ we get
\begin{align*}
\mathcal{F}(y,z)&=y\sum_{k=0}^{q-1}\alpha_k\sum_{j=0}^kC_{2k+1}^{2j+1}y^{2j}z^{2k-2j}\\
&\leq y\max_{0\leq k\leq q-1}\vert\alpha_k\vert 2^{2q}\max_{0\leq j\leq q}\sum_{j=0}^ky^{2j}z^{2k-2j}
\leq Cy\left(1+y^{2q}\right)\left(1+z^{2q}\right).
\end{align*}
Finally since $p>4q$, there exists some constant $C'>0$ such that:
\begin{eqnarray*}
\int_{0}^\infty\mathcal{F}(y,z)u(z)dz&\leq& C y\left(1+y^{2q}\right)\int_{0}^\infty\left(1+z^{2q}\right)u(z)dz\\
&\leq&Cy\left(1+y^{2q}\right)\int_{0}^\infty\frac{1+z^{2q}}{1+z^p}\Big((1+z^p)u(z)\Big)dz\\
&\leq&C' y M\left(1+y^{2q}\right).
\end{eqnarray*}
By integration we obtain
$\displaystyle
\int_0^x(F'\ast u)(y) d y\leq C_2Mx^2(1+x^{2q})$ for all $x\in\mathbb{Re}_+$.
\end{proof}
\begin{lem}\label{lem:inclusion}
There exists $M_0>0$ such that for any $M\geq M_0$, $\mathbb{A}^{\epsilon}(\mathbb{C}_{M})\subset\mathbb{C}_{M}$.
\end{lem}
\begin{proof} By construction $\mathbb{A}^{\epsilon}u$ is a non negative even function which satisfies $\int_{\mathbb{Re}}\mathbb{A}^{\epsilon}u(x)dx=1$. It suffices then to prove that:
$$
\sup_{x\in\mathbb{Re}}\left(1+|x|^p\right)\mathbb{A}^{\epsilon}u(x)\leq M.
$$
By \eqref{def:operateur:A} and according to Lemma \ref{lem:prem:major} we obtain some lower bound for the normalization factor:
\begin{eqnarray*}\lambda_{\epsilon}(u)&=&\int_{-\infty}^{+\infty}\exp\left[-\frac{2}{\epsilon}\left(V(x)+\int_0^x(F'\ast u)(y)dy\right)\right]dx\\
&\geq&\int_{-\infty}^{+\infty}\exp\left[-\frac{2}{\epsilon}\left(V(x)+C_2Mx^2\left(1+x^{2q}\right)\right)\right]dx.
\end{eqnarray*}
According to both (V-3) and (V-7), we know that $V(x)\leq 0$ for all $x\in[-a;a]$. Hence
$$
\lambda_{\epsilon}(u)\geq\int_{-a}^{+a}\exp\left[-\frac{2}{\epsilon}C_2Mx^2(1+a^{2q})\right]dx.
$$
Let us define $\xi(M)=\epsilon^{1/2}(2C_2M(1+a^{2q}))^{-1/2}$ then $\lim_{M\to\infty}\xi(M)=0$. By the change of variable $x:=\xi(M)y$ and Lemma \ref{lem:Annexe1}, the following development holds
\begin{align*}
&\int_{-a}^{+a}\exp\left[-\frac{2}{\epsilon}C_2Mx^2(1+a^{2q})\right]dx=2\xi(M)\int_0^{a/\xi(M)}e^{-x^2}dx\\
&=\xi(M)\left\{\frac{\sqrt{\pi}}{2}-\frac{\xi(M)}{a}\exp\left[-\frac{a^2}{\xi(M)^2}\right]+o\left(\frac{\xi(M)}{a}\exp\left[-\frac{a^2}{\xi(M)^2}\right]\right)\right\}.
\end{align*}
As soon as $M$ is large enough, we have $\lambda_{\epsilon}(u)\geq \sqrt{\pi}\xi(M)\frac{1}{4}=\sqrt{\frac{\pi\epsilon}{32 C_2(1+a^{2q})}}\frac{1}{\sqrt{M}}$.\\
Therefore $\frac{1}{\lambda_{\epsilon}(u)}\leq C(\epsilon)\sqrt{M}$ where $C(\epsilon)$ is a positive constant determined by parameters of the global system and $\epsilon$. By \eqref{def:operateur:A} and the preceding upper bound, we prove that
$$
(1+|x|^p)\mathbb{A}^{\epsilon}u(x)\leq C(\epsilon)\sqrt{M}(1+|x|^p)e^{-\frac{2}{\epsilon}V(x)}\leq C'(\epsilon)\sqrt{M},
$$
where $C'(\epsilon)$ is a positive constant similar to $C(\epsilon)$. In order to conclude, it is sufficient to choose $M\geq C'(\epsilon)^2$: we get immediately $\mathbb{A}^{\epsilon}u\in\mathbb{C}_M$.
\end{proof}
\begin{lem}\label{lem:continu}
$\mathbb{A}^{\epsilon}$ is a continuous operator on $\mathbb{C}_M$ with respect to the uniform norm.
\end{lem}
\begin{proof} We shall find some upper bound for the following expression $\vert\mathbb{A}^\epsilon u-\mathbb{A}^\epsilon v\vert$.\\
{\bf Step 1.} Let $u,v\in\mathbb{C}_{M}$. We define:
\begin{align*}
\Lambda^{\epsilon}(x)=e^{-\frac{2}{\epsilon}V(x) }
\left\{\exp\left[-\frac{2}{\epsilon}\int_0^x(F'\ast u)(y)dy\right]
-\exp\left[-\frac{2}{\epsilon}\int_0^x(F'\ast v)(y)dy\right]\right\}\\
=e^{-\frac{2}{\epsilon}V(x)-\frac{\alpha}{\epsilon}x^2}\left\{\exp\left[-\frac{2}{\epsilon}\int_0^x(F'_0\ast u)(y)dy\right]-\exp\left[-\frac{2}{\epsilon}\int_0^x(F'_0\ast v)(y)dy\right]\right\}\nonumber.
\end{align*}
It is well known that $|e^{-a}-e^{-b}|\leq|a-b|$ for $a,b\geq 0$. In order to apply this inequality we have to prove that $\int_0^x(F'_0\ast v)(y)dy$ and $\int_0^x(F'_0\ast u)(y)dy$ are non negative. By Lemma \ref{lem:prem:major}, for each function $f\in\mathbb{C}_M$ the convolution term  $\int_0^x\left(F'\ast f\right)(y)dy$ is lower bounded by $\frac{\alpha}{2}x^2$. So $\int_0^x\left(F'_0\ast f\right)(y)dy$ is non negative due to the relation: $F'(y)=F'_0(y)+\alpha y$. Hence
\begin{equation}\label{theta}
|\Lambda^{\epsilon}(x)|\leq\frac{2}{\epsilon}e^{-\frac{2}{\epsilon}V(x)-\frac{\alpha}{\epsilon}x^2}\Lambda_0^\epsilon(x),
\end{equation}
with $\Lambda_0^\epsilon$ defined by
\[
\left|\int_0^x(F'_0\ast u)(y)dy-\int_0^x(F'_0\ast v)(y)
dy\right|=\left|\int_0^x\int_{\mathbb{Re}}F'_0(y-z)(u(z)-v(z))dz dy\right|.
\]
Since $u$ and $v$ are elements of $\mathbb{C}_M$, they are even functions and the integral with respect to the variable $z$ becomes
\begin{align}\label{theta_0}
\Lambda_0^\epsilon&=\left\vert\int_0^x\int_0^{\infty}\left(F'_0(z+y)-F'_0(z-y)\right)\left(u(z)-v(z)\right)dz dy\right\vert\nonumber\\
&\le \int_0^x\int_{0}^\infty\left|F'_0(z+y)-F'_0(z-y)\right|\left|u(z)-v(z)\right|dz dy.
\end{align}
The assumption (F-4) gives informations about the increments of the interaction function: there exist two positive constants $C_q$ and $C$ such that
\begin{eqnarray}\label{eq:accroissm}
\left|F'_0(z+y)-F'_0(z-y)\right|&\leq&2|y|C_q\left(1+|z+y|^{2q-2}+|z-y|^{2q-2}\right)\nonumber\\
&\leq&2|y|C_q\left(1+2^{2q-1}|z|^{2q-2}+2^{2q-1}|y|^{2q-2}\right)\nonumber\\
&\leq& C|y|\left(1+|y|^{2q-1}+|z|^{2q-1}\right)\nonumber\\
&\leq& C|y|\left(1+|y|^{2q-1}\right)\left(1+|z|^{2q-1}\right).
\end{eqnarray}
We shall now find some upper bound for $|u(z)-v(z)|$ in \eqref{theta_0}. Since $u,v\in\mathbb{C}_M$ then $u(z)(1+|z|^p)\leq M$ and $v(z)(1+|z|^p)\leq M$, $\forall z\in\mathbb{Re}$. The obvious upper bound $|u(z)-v(z)|(1+|z|^p)\leq 2M$ permits to obtain $\sqrt{|u(z)-v(z)|}\leq\sqrt{\frac{2M}{1+|z|^p}}$. Consequently, for all $z$ of $\mathbb{Re}$, $|u(z)-v(z)|\leq\sqrt{||u-v||_{\infty}}\sqrt{\frac{2M}{1+|z|^p}}$ where $\Vert\cdot\Vert_\infty$ denotes the uniform norm. Using this inequality, \eqref{eq:accroissm} and \eqref{theta_0} in order to estimate $\Lambda_0^\epsilon$, we get
\[
|\Lambda_0^{\epsilon}(x)|\leq C\sqrt{||u-v||_{\infty}}\int_0^x|y|\left(1+|y|^{2q-1}\right)dy\int_{0}^\infty
\sqrt{\frac{2M}{1+|z|^p}}\left(1+z^{2q-1}\right)dz.
\]
Since $p>4q$ the integral with respect to the variable $z$ is finite and can be considered like a constant term. By \eqref{theta} and using the positivity of $\alpha x^2$, we obtain directly the existence of some positive constant $C>0$ such that
$$
|\Lambda^{\epsilon}(x)|\leq C\sqrt{\frac{M}{\epsilon}}\sqrt{||u-v||_{\infty}}x^2\left(1+|x|^{2q-1}\right)e^{-\frac{2}{\epsilon}V(x)}.
$$
According to (V-2), the expression $x^2\left(1+|x|^{2q-1}\right)e^{-\frac{2}{\epsilon}V(x)}$ can be bounded by some constant independent of $\epsilon$. Therefore 
\[
||\Lambda^{\epsilon}||_{\infty}\leq C(M,\epsilon)\sqrt{||u-v||_{\infty}}.
\]
Two results can be deduced: firstly $||\Lambda^{\epsilon}||_{\infty}$ is finite and secondly $||\Lambda^{\epsilon}||_{\infty}$ becomes small as $||u-v||_{\infty}$ decreases towards $0$.\\
{\bf Step 2.} For any $x\in\mathbb{Re}$, we introduce:
\begin{align}\label{eq:omega}
\Omega_{\epsilon}(x)=\frac{1}{\lambda_{\epsilon}(u)\lambda_{\epsilon}(v)}\exp\left[-\frac{2}{\epsilon}\left(\int_0^x(F'\ast v)(y)dy+V(x)\right)\right].
\end{align}
Then the difference $\mathbb{A}^{\epsilon}u(x)-\mathbb{A}^{\epsilon}v(x)$ can be decomposed as follows:
\begin{equation}\label{eq:decomp}
\mathbb{A}^{\epsilon}u(x)-\mathbb{A}^{\epsilon}v(x)=\frac{1}{\lambda_{\epsilon}(u)}\Lambda^{\epsilon}(x)+(\lambda_{\epsilon}(v)-\lambda_{\epsilon}(u))\Omega_{\epsilon}(x).
\end{equation}
Taking the uniform norm, we get 
\begin{equation}\label{eq:decomp:major}
||\mathbb{A}^{\epsilon}u-\mathbb{A}^{\epsilon}v||_{\infty}\leq\frac{1}{\lambda_{\epsilon}(u)}||\Lambda^{\epsilon}||_{\infty}+\left|\lambda_{\epsilon}(v)-\lambda_{\epsilon}(u)\right|\left|\left|\Omega_{\epsilon}\right|\right|_{\infty}.
\end{equation}
We have shown in the proof of Lemma \ref{lem:inclusion} that  $\frac{1}{\lambda_{\epsilon}(u)}\leq C(\epsilon)\sqrt{M}$ and moreover $||\Lambda^{\epsilon}||_{\infty}\leq C(M,\epsilon)\sqrt{||u-v||_{\infty}}$. We deduce that 
\[
\frac{1}{\lambda_{\epsilon}(u)}||\Lambda^{\epsilon}||_{\infty}\leq C'(M,\epsilon)\sqrt{||u-v||_{\infty}}.
\]
It is then sufficient to find a similar inequality for the term $\left|\lambda_{\epsilon}(v)-\lambda_{\epsilon}(u)\right|\left|\left|\Omega_{\epsilon}\right|\right|_{\infty}$ in order to conclude the proof.
\begin{align*}
|\lambda_{\epsilon}(v)-\lambda_{\epsilon}(u)|&=\left|\int_{\mathbb{Re}}\Lambda^{\epsilon}(x)dx\right|\\
&\leq C\sqrt{\frac{M}{\epsilon}}\sqrt{||u-v||_{\infty}}\int_{-\infty}^{+\infty}x^2\left(1+|x|^{2q-1}\right)e^{-\frac{2}{\epsilon}V(x)}dx.
\end{align*}
According to (V-4), the integral with respect to the variable $x$ is finite and does not depend on $M$.
We have immediately 
\[
|\lambda_{\epsilon}(v)-\lambda_{\epsilon}(u)|\leq C(M,\epsilon)\sqrt{||u-v||_{\infty}}.
\]
It remains to estimate $\Omega_\epsilon(x)$. By (V-4) and \eqref{eq:lem:prem:major}, we have
$$\int_0^x\left(F'\ast v\right)(y)dy+V(x)\geq C_4x^4+\Big(\frac{\alpha}{2}-C_2\Big)x^2$$
for all $x$ positive. Furthermore the symmetry property of $V$ and $F$ permits to extend the bound to all $x\in\mathbb{Re}$. The function $\exp\left[-\frac{2}{\epsilon}\left(\int_0^x(F'\ast v)(y)dy+V(x)\right)\right]$ is then bounded by a constant depending on $\epsilon$. Moreover we have already proved that $\frac{1}{\lambda_{\epsilon}(f)}\leq C(\epsilon)\sqrt{M}$ for all elements $f$ of the function space $\mathbb{C}_M$. This bound can therefore be applied to $u$ and $v$. Finally
we obtain the existence of some constant $C(\epsilon)>0$ such that, for all real value $x$, $\left|\Omega_{\epsilon}(x)\right|\leq C(\epsilon)M$. \\
By \eqref{eq:decomp}, we have
$$||\mathbb{A}^{\epsilon}u-\mathbb{A}^{\epsilon}v||_{\infty}\leq C'(M,\epsilon)\sqrt{||u-v||_{\infty}}+C(M,\epsilon)\sqrt{||u-v||_{\infty}}C(\epsilon)M.$$
In other words,
$$||\mathbb{A}^{\epsilon}u-\mathbb{A}^{\epsilon}v||_{\infty}\leq C''(M,\epsilon)\sqrt{||u-v||_{\infty}}$$
what finishes the proof. Here $C$, $C'$ and $C''$ are positive constants.
\end{proof}
We have now all the keys for proving the existence of some symmetric invariant measure. Indeed we have just presented some continuous mapping which stabilizes a convex subset of the Banach space $\mathbb{B}$.
\begin{thm}\label{thm:exist-sym}
There exists a symmetric invariant measure for \eqref{eq:equation_depart}.
\end{thm}
\begin{proof}
Let $M_0$ defined by Lemma \ref{lem:inclusion}. Taking $M\geq M_0$, let us prove that $\overline{\mathbb{A}^{\epsilon}(\mathbb{C}_{M})}$ is a compact set. For this reason we shall estimate the following derivative:
\[
\left(\mathbb{A}^{\epsilon}u\right)'(x)=-\frac{2}{\epsilon}\frac{(F'\ast u)(x)+V'(x)}{\lambda_{\epsilon}(u)}\exp\left[-\frac{2}{\epsilon}\left(\int_0^x(F'\ast u)(y)dy+V(x)\right)\right].
\]
Let us analyze the different elements of this derivative. We have already seen in the proof of Lemma \ref{lem:inclusion} that for any $u\in\mathbb{C}_M$ the normalization factor $\lambda_\epsilon(u)$ satisfies
\begin{equation}\label{et1}
\frac{1}{\lambda_{\epsilon}(u)}\leq C(\epsilon)\sqrt{M}.
\end{equation}
By \eqref{eq:lem:prem:major}, we obtain the bound: $0\leq\int_0^x\left(F'\ast u\right)(y)dy\leq C_2Mx^2\left(1+x^{2q}\right)$.\\
Furthermore by (V-4) and (V-7), we get some estimation of $V$ and its derivative:
\begin{equation}\label{et2}
V(x)\geq C_4x^4-C_2x^2\quad\mbox{and}\ \left|V'(x)\right|\leq C_q\left(1+|x|^{2q}\right)\quad\mbox{for all}\ x\in\mathbb{Re}.
\end{equation}
It remains to find some upper bound for the convolution term: $\left|\left(F'\ast u\right)(x)\right|$ with $x\in\mathbb{Re}_+$. By (F-4) and since $u$ is an even function,
\begin{align*}
\left|\left(F'\ast u\right)(x)\right|=\left|\int_{\mathbb{Re}}F'(x-z)u(z)dz\right|\leq\int_{0}^\infty\Big|F'(x+z)+F'(x-z)\Big|u(z)dz\\
\leq C_q\int_{0}^\infty\Big\{|x+z|\left(1+|x+z|^{2q-2}\right)+|x-z|\left(1+|x-z|^{2q-2}\right)\Big\}u(z)dz.
\end{align*}
 Therefore:
\begin{eqnarray*}
\left|\left(F'\ast u\right)(x)\right|&\leq&\int_{\mathbb{Re}^+}C_q2^{2q-1}\Big\{|x|^{2q-1}+|z||x|^{2q-2}\\
&&+|x|\left(1+|z|^{2q-2}\right)+|z|\left(1+|z|^{2q-2}\right)\Big\}u(z)dz.
\end{eqnarray*}
By definition of $\mathbb{C}_M$, we have $u(z)\leq \frac{M}{1+|z|^p}$ for $p>4q$. Hence the moments of order $1$, $2q-2$ and $2q-1$ are bounded: there exist some constants $C$ and $C'$, independent of the different parameters appearing in the system, such that 
\begin{equation}\label{et3}
\left|\left(F'\ast u\right)(x)\right|\leq C\left(1+|x|+|x|^{2q-2}+|x|^{2q-1}\right)\le C'\left(1+|x|^{2q+1}\right).
\end{equation}
To sum up: using \eqref{et1}, \eqref{et2} and \eqref{et3} we obtain
$$
\left|\left(\mathbb{A}^{\epsilon}u\right)'(x)\right|\leq \frac{2}{\epsilon}C(\epsilon)\sqrt{M}(1+|x|^{2q+1})\exp\left[-\frac{2}{\epsilon}\left(C_4x^4-C_2x^2\right)\right].
$$
Finally we deduce that there exists some constant $C_\epsilon$ such that $\left|\left(\mathbb{A}^{\epsilon}u\right)'(x)\right|\leq C_\epsilon$ for all $x\in\mathbb{Re}$.\\
Let us prove now that $\overline{\mathbb{A}^{\epsilon}\mathbb{C}_M}$ is compact. To this end, we take some sequence of functions $\left(u_n\right)_{n\in\mathbb{N}}$ in $\mathbb{C}_{M}$ and focus our attention to the sequence $\left(\mathbb{A}^{\epsilon}u_n\right)_{n\in\mathbb{N}}$.
According to the definition of $\mathbb{A}^{\epsilon}$, for all $x$ real the set $\left\{\mathbb{A}^{\epsilon}u_n(x),n\in\mathbb{N}\right\}$ is compact. Furthermore the bound of $\left|\left(\mathbb{A}^{\epsilon}u\right)'(x)\right|$ is independent of the variables $x$ and $u\in\mathbb{C}_M$: the equicontinuity condition for the application of Ascoli's theorem is satisfied. Hence, we deduce that there exists some subsequence of $\mathbb{A}^{\epsilon}u_n$ which converges to a limit function $v$ belonging to $\overline{\mathbb{A}^{\epsilon}\mathbb{C}_M}$.\\
By Lemma \ref{lem:inclusion} and Lemma \ref{lem:continu} we can apply Schauder's theorem (Proposition \ref{schauder}) for the operator $\mathbb{A}^{\epsilon}$ on the function space $\mathbb{C}_M$ with $M\ge M_0$. We deduce the existence of some fixed point which is, by construction, a symmetric stationary measure for the diffusion \eqref{eq:equation_depart}.
\end{proof}\mathversion{bold}
\subsection{Example: $F(x)=\frac{\beta}{4}x^4+\frac{\alpha}{2}x^2$}
\mathversion{normal}
We have just shown the existence of a symmetric invariant measure for general self-stabilizing diffusions using fixed point arguments. Now let us study some particular case by a completely different way: the procedure shall be close to that developed in section \ref{sect:set}. Let $V$ be a potential satisfying (V-1)-(V-7).\\
Let $u_{\epsilon}$ be a symmetric invariant measure (Theorem \ref{thm:exist-sym}). We denote by $m_2(\epsilon)$ its second moment. The couple $(m_2(\epsilon),u_\epsilon)$ is solution to some system like \eqref{eq:moyenne}-\eqref{eq:structure:lineaire}. Indeed
\begin{eqnarray*}
F\ast u_{\epsilon}(x)&=&\int_{\mathbb{R}}F(x-z)u_{\epsilon}(z)dz\\
&=&\frac{\alpha}{2}x^2+\frac{\beta}{4}x^4+\frac{3\beta m_2(\epsilon)}{2}x^2+\left(\frac{\alpha}{2}m_2(\epsilon)+\frac{\beta}{4}\int_{\mathbb{R}}z^4u_{\epsilon}(z)dz\right),
\end{eqnarray*}
with $\beta\ge 0$ since $F'$ is a convex function on $\mathbb{Re}_+$.\\
The expression delimited by the brackets is just a constant so we obtain the following system of equations for $m_2(\epsilon)$ and $u_\epsilon$: $m_2(\epsilon)=\int_\mathbb{Re}x^2 \nu(m_2(\epsilon),x)dx$ and $u_\epsilon(x)=\nu(m_2(\epsilon),x)$ where
\begin{eqnarray*}
\nu(m,x)=\frac{\exp\left[-\frac{2}{\epsilon}\left(V(x)+F(x)+\frac{3\beta m}{2}x^2\right)\right]}{\int_0^{\infty}\exp\left[-\frac{2}{\epsilon}\left(V(z)+F(z)+\frac{3\beta m}{2}z^2\right)\right]dz}.
\end{eqnarray*}
Therefore we introduce the function $\chi_{\epsilon}(m)=\int_0^{\infty}x^2\nu(m,x)dx-m$.
By Theorem \ref{thm:exist-sym}, we know that $\chi_{\epsilon}$ admits at least one zero on $\mathbb{Re}_+$. Computing the derivative of $\chi_{\epsilon}$, we prove that the considered function is decreasing:
\[
\chi_{\epsilon}'(m)=-\frac{3\beta}{2}\left\{\int_0^{\infty}x^4\nu(m,x)dx-\left(\int_0^{\infty}x^2\nu(m,x)dx\right)^2\right\}-1<0.
\]
The conclusion is immediate: there is a unique symmetric invariant measure. Obviously this result and the kind of method used to prove it are particular to our simple example. If the degree of the interaction function is strictly larger than 4 then it isn't enough to know the second moment in order to define the invariant measure: we need more moments and the proof of the uniqueness becomes awkward. 
\subsection{Outlying invariant measures}
This section is essentially motivated by the uniqueness question for invariant measures. The existence of some symmetric measure was just proved in Section \ref{sym}. It suffices now to point out asymmetric stationary measures for self-stabilizing diffusions. In the general setting, the interaction function is polynomial: set $F(x)=\sum_{k=1}^n\frac{F^{(2k)}(0)}{(2k)!}x^{2k}$.\\
Let $u$ be the density of some probability measure with respect to the Lebesgue measure and $\mu_1,\cdots,\mu_{2n-1}$ denote its moments of orders $1$ to $2n-1$ respectively. We assume they are finite. Then the difference $D(x):=F\ast u(x)-F\ast u(0)$ satisfies
\begin{align*}
D(x)&=F(x-a)-F(-a)+\sum_{p=1}^{2n-1}\frac{(-1)^p}{p!}(\mu_p-a^p)\sum_{j\geq\frac{1+p}{2}}^n\frac{F^{(2j)}(0)}{(2j-p)!}x^{2j-p}\\
&=F(x-a)-F(a)+\sum_{p=1}^{2n-1}\frac{(-1)^p}{p!}(\mu_p-a^p)\left(F^{(p)}(x)-F^{(p)}(0)\right).
\end{align*}
Hence $D(x)=Z_m(x)-Z_m(0)$ where
\begin{equation}\label{def:zn}
Z_m(x)=F(x-a)+\sum_{p=1}^{2n-1}\frac{(-1)^p}{p!}(m_p-a^p)F^{(p)}(x).
\end{equation}
Since the convolution product can be expressed as a polynomial function which coefficients just depend on the moments of $u$, then the exponential expression of invariant measure \eqref{eq:lem:structure} can be specified. Indeed equation \eqref{eq:lem:structure} can be transformed into some system of equations whose unknown factors are the moments of the measure.
In order to introduce this system, let us define, for all $k\in[1;2n-1]$, the function
\begin{eqnarray}\label{eq:chaq-elem}
\varphi_k^{(\epsilon)}(m_1,\cdots,m_{2n-1})&=&\frac{\int_{\mathbb{Re}}x^k\exp\left[-\frac{2}{\epsilon}\left(V(x)+Z_m(x)-Z_m(0)\right)\right]dx}{\int_{\mathbb{Re}}\exp\left[-\frac{2}{\epsilon}\left(V(x)+Z_m(x)-Z_m(0)\right)\right]dx}\nonumber\\
&=&\frac{\int_{\mathbb{Re}}x^k\exp\left[-\frac{2}{\epsilon}W_m(x)\right]dx}{\int_{\mathbb{Re}}\exp\left[-\frac{2}{\epsilon}W_m(x)\right]dx}
\end{eqnarray}
with the potential $W_m(x)=V(x)+Z_m(x)$. We construct the mapping:
\begin{equation}\label{eq:ensemb}
\Phi^{(\epsilon)}=(\varphi_1^{(\epsilon)} ,\ldots , \varphi_k^{(\epsilon)} , \ldots , \varphi_{2n-1}^{(\epsilon)}).
\end{equation}
\emph{The measure associated to the density function $u$ is invariant if and only if its moments vector $(\mu_1,\cdots,\mu_{2n-1})$ is a fixed point of the map $\Phi^{(\epsilon)}$.}

\noindent We are going to show the existence of an asymmetric invariant measure defined by $2n-1$ parameters close to $a,\cdots,a^{2n-1}$ respectively, in other words the outlying measure is close to the Dirac mass in the point $a$.
More precisely, we shall prove that there exists a parallelepiped stable by $\Phi^{(\epsilon)}$, which converges to the point $(a,a^2,\cdots,a^{2n-1})$ as $\epsilon$ tends to $0$. 
As in the linear case, we shall proceed by applying the mean value theorem in order to obtain asymptotic developments in the small noise limit.
\begin{thm}\label{thm:exist:asym} Let $(\eta_\epsilon)_\epsilon$ some sequence satisfying $\displaystyle\lim_{\epsilon\to 0}\eta_\epsilon=0$ and $\displaystyle\lim_{\epsilon\to 0}\epsilon/\eta_\epsilon=0$. Under the condition
\begin{equation}\label{eq:thm:exist:asym}
\sum_{p=0}^{2n-2}\frac{\left|F^{(p+2)}(a)\right|}{p!}\ a^{p}<\alpha+V''(a),
\end{equation}
 for any $\rho>0$, there are at least two outlying measures $u_\epsilon^+$ and $u_\epsilon^-$ satisfying, for $\epsilon$ small enough
\begin{equation}
 \label{eq:thm:exist:asym2}
\left\vert \int_\mathbb{Re} x^ku_\epsilon^\pm (x)dx-(\pm a)^k\right\vert\le \rho\ \eta_\epsilon.
\end{equation}

\end{thm}
\begin{proof}Let $\lambda>0$.
Let us define the parallelepiped 
\[
C(\epsilon)=\prod_{p=1}^{2n-1}[a^p-pa^{p-1}\lambda\eta_\epsilon,a^p+pa^{p-1}\lambda\eta_\epsilon].
\]
Let $m$ be an element of $C(\epsilon)$ then there exist some coordinates $(r_p)_{1\le p\le 2n-1}$ which determine $m$ through the equations $m_p=a^p+r_p\,\eta_\epsilon$. By \eqref{def:zn} and \eqref{eq:chaq-elem}, we get
\[
\varphi_k^{(\epsilon)}(m)=\frac{\int_{\mathbb{Re}}x^ke^{-\frac{2}{\epsilon}(V(x)+F(x-a))}\exp\Big[-\frac{2\eta_\epsilon}{\epsilon}\sum_{p=1}^{2n-1}\frac{(-1)^pr_p}{p!}F^{(p)}(x)\Big]dx}{\int_{\mathbb{Re}}e^{-\frac{2}{\epsilon}(V(x)+F(x-a))}\exp\Big[-\frac{2\eta_\epsilon}{\epsilon}\sum_{p=1}^{2n-1}\frac{(-1)^pr_p}{p!}F^{(p)}(x)\Big]dx}.
\]
We apply Lemma \ref{lem:Annexe3} and Remark \ref{lem:rem1} to the functions $U(x)=V(x)+F(x-a)$, $f(x)=x^k$, $\mu_p=r_p$ and $G_p(x)=\frac{(-1)^p}{p!}F^{(p)}(x)$. We obtain:
\[
\varphi_k^{(\epsilon)}(m)=a^k-\eta_\epsilon\frac{k a^{k-1}}{\alpha+V''(a)}\sum_{p=1}^{2n-1}\frac{(-1)^pr_p}{p!} F^{(p+1)}(a)+o(\eta_\epsilon),
\]
uniformly with respect to the coordinates $(r_p)_p$.
By definition of the parallelepiped $C(\epsilon)$ the coordinates satisfy $\left|r_p\right|\leq pa^{p-1}\lambda$. Therefore, under condition \eqref{eq:thm:exist:asym},
\begin{eqnarray*}
\left|\varphi_k^{(\epsilon)}(m)-a^k\right|&\leq&\eta_\epsilon\lambda\frac{ka^{k-1}}{\alpha+V''(a)}\sum_{p=1}^{2n-1}\frac{\left|F^{(p+1)}(a)\right|}{p!}pa^{p-1}+o(\eta_\epsilon)\\
&<&\eta_\epsilon ka^{k-1}\lambda+o(\eta_\epsilon).
\end{eqnarray*}
Since this estimate is uniform with respect to the coordinates, 
 as soon as $\epsilon$ is small enough, we have
$|\varphi_k^{(\epsilon)}(m)-a^k|<ka^{k-1}\lambda\eta_\epsilon$,
that means that $\Phi^{(\epsilon)}(m)\in C(\epsilon)$.\\
Let us note that $C(\epsilon)$ is a convex, closed and bounded subset of $\mathbb{Re}^{2n-1}$. Since the space dimension is finite, the continuity of $\Phi^{(\epsilon)}$ implies that the closure of the parallelepiped's image is a compact set.\\
We can apply Schauder's Theorem (Proposition \ref{schauder}) and obtain that there exists some fixed point in the compact. In other words there exists $m\in C(\epsilon)$ such that the measure associated to the density
\begin{equation}\label{assoc}
u_{\epsilon,m}(x)=\frac{\exp\left[-\frac{2}{\epsilon}W_m(x)\right]}{\int_{\mathbb{R}}\exp\left[-\frac{2}{\epsilon}W_m(z)\right]dz}
\end{equation}
is invariant. In a similar way, the measure defined by $m^-$ is also invariant; here  $m^-(k)=(-1)^km_k$. To conclude: we have at least two outlying measures, one around $a$ and the second one around $-a$.
\end{proof}
We can not prove at this stage the uniqueness of the couple of outlying invariant measures (this question shall be explored in a subsequent work). We can effectively imagine that other outlying measures could exist around $a$, around $-a$ or even around other areas. Nevertheless we can develop a sharper description of one particular outlying measure: the measure close to $\delta_a$ where $\delta$ represents the Dirac measure. To do this it suffices to estimate its different moments, that requires the following preliminary result. 
\begin{lem}\label{lem:sol:syst}
There exists a unique solution $(\tau_1^0,\cdots,\tau_{2n-1}^0)$ to the following Cramer's system
\begin{equation}\label{eq:cramer}
\sum_{p=1}^{2n-1}\frac{(-1)^p}{p!}\ F^{(p+1)}(a)\tau_p+\frac{\alpha+V''(a)}{ka^{k-1}}\ \tau_k=\frac{V^{(3)}(a)}{4(\alpha+V''(a))}-\frac{k-1}{4a},
\end{equation}
for $1\le k\le 2n-1$. This solution is given by
\begin{equation}\label{eq:sol:cramer}
\tau_k^0=k a^{k-1}\frac{a V^{(3)}(a)-(k-1)V''(a)}{4aV''(a)\left(\alpha+V''(a)\right)},\quad 1\le k\le 2n-1.
\end{equation}
\end{lem}
\begin{proof} Let us denote by $I_{2n-1}$ the unit matrix of dimension $2n-1$ and for $A\in\mathbb{Re}^{2n-1}$, $A^{\textsc{T}}$ represents the transpose of the vector $A$. Moreover we adopt the following notation $(x_k)_{1\le k\le 2n-1}=(x_1,\ldots,x_{2n-1})$. 
The system \eqref{eq:cramer} can be written in this way: we define $\mathcal{T}=(\tau_k)_{1\le k\le 2n-1}^{\textsc{T}}$ then
\begin{align*}
\Big[(\alpha+V''(a))I_{2n-1}+C_1C_2^{\textsc{T}}\Big]\mathcal{T}=\Big(k a^{k-1}\Big( \frac{V^{(3)}(a)}{4\alpha+4V''(a))}-\frac{k-1}{4a}\Big)\Big)_{{\scriptstyle 1\le k\le 2n-1}}^{\textsc{T}}
\end{align*}
with the vectors $C_1^{\textsc{T}}=(ka^{k-1})_{1\le k\le 2n-1}$ and $C_2^{\textsc{T}}=\Big(\frac{(-1)^k}{k!}F^{(k+1)}(a)\Big)_{1\le k\le 2n-1}$.
We define therefore
\begin{equation}\label{def:A}
A=(\alpha+V''(a))I_{2n-1}+C_1C_2^{\textsc{T}}.
\end{equation}
Let us note that $C_1C_2^{\textsc{T}}C_1C_2^{\textsc{T}}=(C_2^{\textsc{T}}C_1)C_1C_2^{\textsc{T}}$ and
\begin{eqnarray*}
C_2^{\textsc{T}}C_1&=&\sum_{p=1}^{2n-1}\frac{(-1)^p}{p!}F^{(p+1)}(a)pa^{p-1}=-\sum_{p=0}^{2n-2}\frac{(-1)^p}{p!}F^{(p+2)}(a)a^p=-F''(0).
\end{eqnarray*} 
Since $F''(0)=\alpha$, we obtain
\begin{align*}
A^2&=(\alpha+V''(a))^2I_{2n-1}+\left(2(\alpha+V''(a))+C_2^{\textsc{T}}C_1\right)C_1C_2^{\textsc{T}}\\
&=(\alpha+V''(a))^2I_{2n-1}+\Big(2(\alpha+V''(a))-F''(0)\Big)C_1C_2^{\textsc{T}}\\
&=(\alpha+V''(a))^2I_{2n-1}+\Big(\alpha+2V''(a)\Big)C_1C_2^{\textsc{T}}\\
&=(\alpha+2V''(a))A-V''(a)\left(\alpha+V''(a)\right)I_{2n-1},
\end{align*}
 We deduce that $A$ is invertible, that is \eqref{eq:cramer} is a Cramer's system, and using \eqref{def:A} we get explicitly the inverse:
 \begin{eqnarray*}
A^{-1}&=&\frac{1}{V''(a)(\alpha+V''(a))}\Big((\alpha+2V''(a))I_{2n-1}-A\Big)\\
&=&\frac{1}{V''(a)(\alpha+V''(a))}\left(V''(a)I_{2n-1}-C_1C_2^{\textsc{T}}
\right).
\end{eqnarray*}
Therefore the Cramer's system \eqref{eq:cramer} admits a unique solution given by
\begin{align*}
\tau_k^0&=\frac{1}{V''(a)(\alpha+V''(a))}\left\{V''(a)ka^{k-1}\frac{aV^{(3)}(a)-(k-1)(\alpha+V''(a))}{4a\left(\alpha+V''(a)\right)}\right.\\
&-\left.ka^{k-1}\sum_{p=1}^{2n-1}\frac{(-1)^p}{p!}F^{(p+1)}(a)pa^{p-1}\frac{aV^{(3)}(a)-(p-1)(\alpha+V''(a))}{4a\left(\alpha+V''(a)\right)}\right\}\\
&=\frac{ka^{k-1}}{4aV''(a)(\alpha+V''(a))^2}\left\{aV^{(3)}(a)\left[V''(a)-\sum_{p=1}^{2n-1}\frac{(-1)^pa^{p-1}}{(p-1)!}F^{(p+1)}(a)\right]\right.\\
&-\left.(\alpha+V''(a))\left[(k-1)V''(a)-\sum_{p=2}^{2n-1}\frac{(-1)^p}{(p-2)!}F^{(p+1)}(a)a^{p-1}\right]\right\}\\
&=ka^{k-1}\frac{aV^{(3)}(a)-(k-1)V''(a)}{4aV''(a)(\alpha+V''(a))}.
\end{align*}
Indeed, we use 
$$\sum_{p=1}^{2n-1}\frac{(-1)^p}{(p-1)!}F^{(p+1)}(a)a^{p-1}=-F''(0)=-\alpha,
$$
and
$$\sum_{p=2}^{2n-1}\frac{(-1)^p}{(p-2)!}F^{(p+1)}(a)a^{p-1}=aF^{(3)}(0)=0.$$
\end{proof}
Theorem \ref{thm:exist:asym} points out the existence of two outlying measures, one concentrated around $a$ and an other around $-a$. According to Lemma \ref{lem:sol:syst} we get some sharper upper bound for the distance between $\delta_a$  and some asymmetric invariant measure.
\begin{thm}\label{thm:proche}
Under the condition \eqref{eq:thm:exist:asym}, for any $\delta>0$, there exists $\epsilon_0$ such that $\Phi^{(\epsilon)}$ admits two fixed points $m^{\pm}$ with
\begin{eqnarray}\label{eq:thm:proche}
\left|m_k^{\pm}(\epsilon)-\Big((\pm1)^ka^k-(\pm1)^k\tau_k^0\epsilon\Big)\right|\leq\delta\left|\tau_k^0\right|\epsilon,\quad 1\le k\le 2n-1,\ \epsilon\le\epsilon_0.
\end{eqnarray}
\end{thm}
\begin{proof} It is similar to the proof of Theorem \ref{thm:exist:asym}.\\
Let $\delta>0$ and $C(\epsilon)=\prod_{p=1}^{2n-1}[a^p-(\tau_p^0+pa^{p-1}\delta)\epsilon,a^p-(\tau_p^0-pa^{p-1}\delta)\epsilon]$. We choose an element $m$ in the parallelepiped $C(\epsilon)$. For all $1\le p\le 2n-1$, there exists a coordinate $\delta_p\in[-\delta;\delta]$ such that $m_p=a^p-(\tau_p^0+pa^{p-1}\delta_p)\epsilon$. By \eqref{def:zn} and  \eqref{eq:chaq-elem}, we obtain
\[
\varphi_k^{(\epsilon)}(m)=\frac{\int_{\mathbb{Re}}x^k\exp\left[2\sum_{p=1}^{2n-1}\frac{(-1)^p}{p!}(\tau_p^0+pa^{p-1}\delta_p)F^{(p)}(x)\right]e^{-\frac{2}{\epsilon}(V(x)+F(x-a))}dx}{\int_{\mathbb{Re}}\exp\left[2\sum_{p=1}^{2n-1}\frac{(-1)^p}{p!}(\tau_p^0+pa^{p-1}\delta_p)F^{(p)}(x)\right]e^{-\frac{2}{\epsilon}(V(x)+F(x-a))}dx}
\]
We apply Lemma \ref{lem:Annexe2} and Remark \ref{lem:rem1} with the following functions: $U(x)=V(x)+F(x-a)$, $\mu_p=\tau_{p}^0+pa^{p-1}\delta_p$, $G=0$ and $f_{p}(x)=2\frac{(-1)^p}{p!}F^{(p)}(x)$. Hence 
\begin{align*}
\varphi_k^{(\epsilon)}(m)&=a^k
-\frac{ka^{k-2}}{4(\alpha+V''(a))^2}\Big[aV^{(3)}(a)-(\alpha+V''(a))\Big((k-1)\\
&\quad+4a\sum_{p=1}^{2n-1}\frac{(-1)^p}{p!}(\tau_p^0+pa^{p-1}\delta_p)F^{(p+1)}(a)\Big)\Big]\epsilon+o(\epsilon)\\
&=a^k-\frac{1}{\alpha+V''(a)}\Big[\frac{ka^{k-1}V^{(3)}(a)}{4(\alpha+V''(a))}-k a^{k-1}\sum_{p=1}^{2n-1}\frac{(-1)^p\tau_p^0}{p!}F^{(p+1)}(a)\\
&-\frac{k(k-1)a^{k-2}}{4}-ka^{k-1}\sum_{p=1}^{2n-1}\frac{(-1)^p\delta_p a^{p-1}}{(p-1)!}F^{(p+1)}(a)\Big]\epsilon+o(\epsilon).
\end{align*}
This estimate is uniform with respect to the variables $(\delta_p)_p$.\\
We denote by $d_k^\epsilon$ the difference $|\varphi_k^{(\epsilon)}(m)-a^k+\tau_k^0\epsilon|$. We compute this expression:
\begin{align*}
d_k^\epsilon &\leq\left|\frac{ka^{k-1}V^{(3)}(a)}{4(\alpha+V''(a))^2}-\frac{ka^{k-1}}{\alpha+V''(a)}\sum_{p=1}^{2n-1}\frac{(-1)^p}{p!}\tau_p^0F^{(p+1)}(a)\right.\\
&-\frac{k(k-1)a^{k-2}}{4(\alpha+V''(a))}-\tau_k^0
-\left.\frac{ka^{k-1}}{\alpha+V''(a)}\sum_{p=1}^{2n-1}\frac{(-1)^p}{(p-1)!}\delta_pF^{(p+1)}(a)a^{p-1}\right|\epsilon+o(\epsilon).
\end{align*}
According to the Lemma \ref{lem:sol:syst} and using the condition \eqref{eq:thm:exist:asym}, we obtain, for $\epsilon$ small enough,
\begin{align*}
\left|\varphi_k^{(\epsilon)}(m)-a^k+\tau_k^0\epsilon\right|&\leq\frac{ka^{k-1}}{\alpha+V''(a)}\sum_{p=1}^{2n-1}\frac{a^{p-1}}{(p-1)!}|\delta_p||F^{(p+1)}(a)|\epsilon+o(\epsilon)\\
&\leq\delta\frac{ka^{k-1}}{\alpha+V''(a)}\sum_{p=0}^{2n-2}\frac{1}{p!}|F^{(p+2)}(a)|a^p\epsilon+o(\epsilon)<ka^{k-1}\delta\epsilon.
\end{align*}
In other words, $\Phi^{(\epsilon)}(m)\in C(\epsilon)$ in the small noise limit. The application of Schauder's Theorem (Proposition \ref{schauder}) permits to prove the existence of some fixed point in the compact. Therefore there exists $m\in C(\epsilon)$ such that the associated measure $u_{\epsilon,m}(x)$ defined by \eqref{assoc} is invariant. In the same way, the measure defined by $m^-$ is invariant with $m^-(k)=(-1)^km_k$. Finally the continuous map $\Phi^{(\epsilon)}$ admits two fixed points $m^{\pm}(\epsilon)$ satisfying \eqref{eq:thm:proche}.
\end{proof}
\begin{rem} {\bf 1.} In the particular case: $F^{(p)}(a)\geq0$ for all $p\in\mathbb{N}$, the condition for the existence of outlying measures becomes
\(
V''(a)>F_0''(2a)\)
where $F_0$ is defined by $F(x)=\frac{\alpha}{2}x^2+F_0(x)$.\\
{\bf 2.} In the linear interaction case: $F(x)=\frac{\alpha}{2}x^2$, \eqref{eq:thm:exist:asym} is equivalent to the simple condition 
$V''(a)>0$ which is in fact always satisfied according to (V-3). In other words we obtain the existence result presented in the linear interaction case.
\end{rem}
\appendix\section{Annexe}
We shall present here some useful asymptotic results which are close to the classical Laplace's method.
\begin{lem}\label{lem:Annexe1} Let $M>0$. Let us assume that $U$ is $\mathcal{C}^2([M,\infty[)$-continuous, $U(x)\neq 0$ and $U''(x)>0$ for all $x\in[M,\infty[$ and $\lim_{x\to\infty}\frac{U''(x)}{(U'(x))^2}=0$.
If $x\to e^{-U(x)}$ is integrable on $\mathbb{Re}$ then for any $m\in\mathbb{Re}$:
\begin{equation}\label{eq:lem:Annexe1}
\int_x^{+\infty}e^{-U(t)}dt\approx\frac{e^{-U(x)}}{U'(x)}\quad\mbox{and}\quad\int_m^xe^{U(t)}dt\approx\frac{e^{U(x)}}{U'(x)}\ \mbox{as}\ x\to\infty.
\end{equation}\end{lem}
\begin{proof} Since $x\to e^{-U(x)}$ is integrable and since these properties are satisfied: $U(x)\neq 0$ and $U''(x)>0$ for $x\ge M$, we know that $\lim_{x\to\infty}U(x)=+\infty$. Furthermore there exists some $M_0>M$ such that $U'(x)>0$ for $x\ge M_0$. Hence for $t\ge M_0$ we obtain
\[
e^{-U(t)}=\left(-\frac{e^{-U(t)}}{U'(t)}\right)'-\frac{U''(t)}{(U'(t))^2}e^{-U(t)}.
\]
Therefore 
\[
I:=\int_x^{\infty}e^{-U(t)}dt=\frac{e^{-U(x)}}{U'(x)}-\int_x^{\infty}\frac{U''(t)}{(U'(t))^2}e^{-U(t)}dt,\quad x\ge M_0.
\]
Using the assumptions of the statement we have $\int_x^{\infty}\frac{U''(t)}{(U'(t))^2}e^{-U(t)}dt\geq 0$. Hence $I\leq e^{-U(x)}U'(x)^{-1}$. Moreover $\lim_{x\longrightarrow\pm\infty}\frac{U''(x)}{(U'(x))^2}=0$. As a consequence for any $\delta>0$, there exists $M_1(\delta)>M_0$ such that $(1+\delta)I\geq\frac{e^{-U(x)}}{U'(x)}$. The estimation of $I$ can be deduced easily. The second equivalence can be obtained by similar arguments. 
\end{proof}
\begin{lem}\label{lem:Annexe:inter1}
Set $\epsilon>0$. Let $U$ and $G$ two $\mathcal{C}^{\infty}(\mathbb{Re})$-continuous functions. We define $U_\mu=U+\mu G$ for $\mu$ belonging to some compact interval $\mathcal{I}$ of $\mathbb{Re}$.
Let us introduce some interval $[a,b]$ satisfying: $U'_\mu(a)\neq 0$, $U'_\mu(b)\neq 0$ and $U_\mu(x)$ admits some unique global minimum on the interval $[a,b]$ reached at $x_\mu\in]a,b[$ for all $\mu\in \mathcal{I}$. We assume that there exists some exponent $k_0$ independent of $\mu\in \mathcal{I}$ such that $2k_0=\min_{r\in\mathbb{N}^*}\left\{U^{(r)}_\mu(x_\mu)\neq 0\right\}$. Then taking the limit $\epsilon\to 0$ we get
\begin{equation}\label{eq:lem:Annexe:inter1}
I_0:=\int_a^b e^{-\frac{U_\mu(t)}{\epsilon}}dt= \frac{1}{k_0}\left(\frac{\epsilon(2k_0)!}{U_\mu^{2k_0}(x_\mu)}\right)^{\frac{1}{2k_0}}\Gamma\left(\frac{1}{2k_0}\right)e^{-\frac{U_\mu(x_\mu)}{\epsilon}}(1+o_{\mathcal{I}}(1)),
\end{equation}
where $\Gamma$ represents the Euler function and $o_{\mathcal{I}}(1)$ converges towards $0$ uniformly with respect to $\mu\in \mathcal{I}$. 
\end{lem}
\begin{proof}
We define $\eta_\mu=\frac{U_\mu^{(2k_0)}(x_\mu)}{(2k_0)!}$. Let us note that $\eta_\mu$ depends continuously on $\mu$. Since $U_\mu$ is regular and admits some unique global minimum for $x=x_\mu$, there exists $\tau_0>0$ independent of the parameter $\mu$ such that $\tau_0<\min\left\{x_\mu-a;b-x_\mu\right\}$ for all $\mu\in\mathcal{I}$ and such that the minimum on the interval $[a;x_\mu-\tau]\bigcup[x_\mu+\tau;b]$ denoted by  $\underline{U}_{\mu}(\tau)$ is reached on the boundary  $\left\{x_\mu-\tau;x_\mu+\tau\right\}$ for all $\tau<\tau_0$. Consequently
\[
\int_{a}^{x_\mu-\tau}\exp\left[-\frac{U_\mu(t)}{\epsilon}\right]dt+\int_{x_\mu+\tau}^{b}\exp\left[-\frac{U_\mu(t)}{\epsilon}\right]dt\leq (b-a)\exp\left[-\frac{\underline{U}_{\mu}(\tau)}{\epsilon}\right].
\]	
Defining $I_\tau=\int_{x_\mu-\tau}^{x_\mu+\tau}\exp\left[-\frac{U_\mu(t)}{\epsilon}\right]dt$, we obtain the following bound:	
\begin{equation}\label{eq:locali}
\left\vert I_0-I_\tau\right\vert\le (b-a)\exp\left[-\frac{\underline{U}_\mu(\tau)}{\epsilon}\right].
\end{equation}
Let us first estimate $I_\tau$. By the mean value theorem, there exists some constant $C>0$ independent of $\mu\in \mathcal{I}$ such that, in a neighborhood of $x_\mu$, the following bound is satisfied: $\left|U_\mu(t)-U_\mu(x_\mu)-\eta_\mu(t-x_\mu)^{2k_0}\right|\leq C|t-x_\mu|^{2k_0+1}$. Hence
\begin{equation}\label{eq:11}
J_1\exp\Big[-\frac{C\tau^{2k_0+1}}{\epsilon}\Big]\le\frac{I_{\tau}}{2}\, \exp\frac{U_\mu(x_\mu)}{\epsilon}\le J_1\exp\Big[\frac{C\tau^{2k_0+1}}{\epsilon}\Big],
\end{equation}
where 
\[
J_\tau=\int_0^{\tau}\exp\left[\frac{1}{\epsilon}\eta_\mu t^{2k_0}\right]dt=\left(\frac{\epsilon}{\eta_\mu}\right)^{\frac{1}{2k_0}}\frac{1}{2k_0}\int_0^{\tau^{2k_0}\frac{\eta_\mu}{\epsilon}}t^{\frac{1}{2k_0}-1}e^{-t}dt,
\]
by the change of variable $t:=\left(\frac{\epsilon}{\eta_\mu}\right)^{\frac{1}{2k_0}}(t')^{\frac{1}{2k_0}}$. A simple integration leads to
\begin{equation}\label{encore}
-\tau^{1-2k_0}\left(\frac{\eta_\mu}{\epsilon}\right)^{\frac{1}{2k_0}-1}e^{-\tau^{2k_0}\frac{\eta_\mu}{\epsilon}}\leq\int_0^{\tau^{2k_0}\frac{\eta_\mu}{\epsilon}}t^{\frac{1}{2k_0}-1}e^{-t}dt-\Gamma\left(\frac{1}{2k_0}\right)\le 0.
\end{equation}
In order to conclude we choose a particular value for the variable $\tau$ namely  $\tau=\exp\left[\frac{\log\left(\epsilon\right)}{2k_0+\frac{1}{2}}\right]$. Then we get: for $C\in\mathbb{Re}$, $l>0$,
\[
\lim_{\epsilon\to 0}e^{C\frac{\tau^{2k_0+1}}{\epsilon}}=1,\quad \lim_{\epsilon\to 0}e^{-\eta_\mu\frac{\tau^{2k_0}}{\epsilon}}\frac{\tau^{1-2k_0}}{\epsilon^{\frac{1}{2k_0}-1}}=\lim_{\epsilon\to 0}\epsilon^{-l}e^{\frac{U_\mu(x_\mu)-\underline{U}_\mu(\tau)}{\epsilon}}=0.
\]
These convergences are uniform with respect to the parameter $\mu$.
Applying these asymptotic results to \eqref{eq:locali}, \eqref{eq:11} and \eqref{encore} permits to prove the statement of the lemma.
\end{proof}

\begin{lem}\label{lem:Annexe:inter2} Let $U$ and $G$ be two $\mathcal{C}^{\infty}([a,b])$-functions. We define $U_\mu=U+\mu G$ for $\mu$ belonging to some compact interval $\mathcal{I}$ of $\mathbb{Re}$. We assume that $U_\mu$ admits a unique global minimum on the interval $]a;b[$ reached at $x=x_\mu$, with $U''_\mu(x_\mu)>0$. Let $f_m$ be a $\mathcal{C}^3$-continuous function for any parameter value $m$ belonging to some compact set $\mathcal{M}$. We assume that there exists some constant $\lambda$ such that $|f_m^{(i)}(x)|\le \lambda$ for all $m\in\mathcal{M}$, $x\in[a,b]$ and $0\le i\le 3$. Then the following asymptotic result holds:
\begin{equation}\label{eq:lem:inter2-1}
\int_{a}^bf_m(t)e^{\frac{-2U_\mu(t)}{\epsilon}}dt=
\sqrt{\frac{\pi\epsilon}{\mathcal{U}_2}}\ e^{-\frac{2U_\mu(x_\mu)}{\epsilon}}\Big\{f_m(x_\mu)+\gamma_0(\mu)\epsilon+o_{\mathcal{I}\mathcal{M}}(\epsilon)\Big\}
\end{equation}
with
\begin{equation}\label{eq:lem:inter2-2}
\gamma_0(\mu)=f_m(x_\mu)\left(\frac{5\ \mathcal{U}_3^2}{48\ \mathcal{U}_2^3}-\frac{\mathcal{U}_4}{16\ \mathcal{U}_2^2}\right)-f'_m(x_\mu)\frac{\mathcal{U}_3}{4\ \mathcal{U}_2^2}+\frac{f''_m(x_\mu)}{4\ \mathcal{U}_2}.
\end{equation}
Here $\mathcal{U}_k=U_\mu^{(k)}(x_\mu)$ and $o_{\mathcal{IM}}(\epsilon)/\epsilon$ converges to $0$ as $\epsilon$ becomes small uniformly with respect to the parameters $m$ and $\mu$.
\end{lem}
\begin{proof}
First we split the integral into two parts:
\[
I=\int_{x_\mu-\rho}^{x_\mu+\rho}f_m(t)e^{\frac{-2U_\mu(t)}{\epsilon}}dt+\int_{[x_\mu-\rho;x_\mu+\rho]^c\bigcap[a;b]}f_m(t)e^{\frac{-2U_\mu(t)}{\epsilon}}dt=I_1+I_2
\]
with some arbitrary $\rho>0$ which should be specified in the following. \\
{\bf Step 1.} We shall prove that the second integral is negligible as $\frac{\rho^2}{\epsilon}\to\infty$ that means that
 $I_2=o_{\mathcal{IM}}\{\epsilon^{3/2}e^{-\frac{2U(x_\mu)}{\epsilon}}\}$. We get
\begin{eqnarray}\label{eq:ret1}
I_2\leq\left(b-a\right)\sup_{z\in[a,b]}\vert f_m(z)\vert \exp\left[-2\frac{\inf_{z\in[x_\mu-\rho;x_\mu+\rho]^c}U_\mu(z)}{\epsilon}\right]
\end{eqnarray}
Since the global minimum of $U_\mu$ is unique and due to the regularity of $U_\mu$ with respect to the parameter $\mu$, we deduce that the minimum of the function on the interval $[x_\mu-\rho;x_\mu+\rho]^c\bigcap[a;b]$ is reached on the boundary provided that $\rho$ is small enough. The development $U_\mu(x_\mu\pm\rho)=U_\mu(x_\mu)+\frac{1}{2}U''_\mu(x_\mu)\rho^2+o_{\mathcal{I}}(\rho^2)$ implies, as already claimed that $I_2=o_{\mathcal{IM}}\left\{\epsilon^{3/2}e^{-\frac{2U(x_\mu)}{\epsilon}}\right\}$ as $\rho^2/\epsilon\to \infty$.\\
{\bf Step 2.} Let us focus our attention to the integral on the domain $[x_\mu-\rho;x_\mu+\rho]$. The function $f_m$ can be developed in the neighborhood of $x_\mu$:
\[
f_m(x)=f_m(x_\mu)+f'_m(x_\mu)(x-x_\mu)+\frac{1}{2}f''_m(x_\mu)(x-x_\mu)^2+\frac{1}{6}f^{(3)}_m(w_{m,\mu}(x))(x-x_\mu)^3
\]
with the value $w_{m,\mu}(x)$ between $x_\mu$ and $x$. Taking into account these different terms, the integral $I_1$ can be split into $4$ different integrals respectively $\tilde{I}_0$,...,$\tilde{I}_3$. For each integral we shall analyze the asymptotic behavior.\\
{\bf Step 2.1.} Asymptotic behavior of $\tilde{I}_3$. By definition $w_{m,\mu}(t)\in[x_\mu-\rho;x_\mu+\rho]$ when $t\in[x_\mu-\rho;x_\mu+\rho]$. Moreover, by assumption $|f_m^{(3)}(w_{m,\mu}(t))|$ is upper bounded by some constant $\lambda>0$ independent of $m$ and $\mu$. By Lemma \ref{lem:Annexe:inter1} applied to $2U_\mu$ , for $\rho<1$ and $\epsilon$ small, we obtain the existence of some constant $C>0$, independent of the parameters $m$ and $\mu$, such that
\[
|\tilde{I}_3|\leq \frac{\lambda}{6}\rho^3\int_{(x_\mu-1)\vee a}^{(x_\mu+1)\wedge b}e^{-\frac{2U_\mu(t)}{\epsilon}}dt\le C\sqrt{\pi}\rho^3\sqrt{\frac{\epsilon}{U''_\mu(x_\mu)}}e^{-\frac{2U_\mu(x_\mu)}{\epsilon}}.
\]
Hence, if $\rho^3=o(\epsilon)$ then the following asymptotic result holds 
\begin{equation}\label{eq:step21}
\tilde{I}_3=o_{\mathcal{IM}}\left\{\epsilon^{\frac{3}{2}}e^{-\frac{2U(x_0)}{\epsilon}}\right\}.
\end{equation}
{\bf Step 2.2.} Asymptotic behavior of $\tilde{I}_2$. Using the $\mathcal{C}^3$-regularity of $U_\mu$ that is $U_\mu(t)=U_\mu(x_\mu)+\frac{1}{2}U''_\mu(x_\mu)(t-x_\mu)^2+\frac{1}{6}U_\mu^{(3)}(y_\mu(t))(t-x_\mu)^3$ with $y_\mu(t)$ belonging to $[x_\mu-\rho;x_\mu+\rho]$, we get
	\[
\tilde{I}_2=\frac{f''_m(x_\mu)}{2}e^{-\frac{2U_\mu(x_\mu)}{\epsilon}}\int_{x_\mu-\rho}^{x_\mu+\rho}(t-x_\mu)^2 e^{-\frac{U''_\mu(x_\mu)}{\epsilon}(t-x_\mu)^2-\frac{U_\mu^{(3)}(y_\mu(t))}{3\epsilon}(t-x_\mu)^3}dt.
		\]
Since $y_\mu(t)$ belongs to some compact set, the third derivative $U_\mu^{(3)}(y_\mu(t))$ is bounded by some constant independent of $\mu$. Applying the following change of variable $u=(t-x_\mu)^2 U''_\mu(x_\mu)/\epsilon$ yields 
\[		
J_2 e^{-C\frac{\rho^3}{\epsilon}}\left(\frac{\epsilon}{U''_\mu(x_\mu)}\right)^{\frac{3}{2}}\le\frac{2\tilde{I}_2\ e^{\frac{2U_\mu(x_\mu)}{\varepsilon}}}{f''_m(x_\mu)}\le J_2 e^{C\frac{\rho^3}{\epsilon}}\left(\frac{\epsilon}{U''_\mu(x_\mu)}\right)^{\frac{3}{2}},
\]
with $J_2=\int_{0}^{U''_\mu(x_\mu)\frac{\rho^2}{\epsilon}}\sqrt{u}e^{-u}du$. If $\frac{\rho^3}{\epsilon}\to 0$ and $\frac{\rho^2}{\epsilon}\to\infty$ then
\begin{equation}\label{eq:step22}	\tilde{I}_2=\sqrt{\pi}\frac{f''_m(x_\mu)}{4}e^{-\frac{2U_\mu(x_\mu)}{\epsilon}}\left(\frac{\epsilon}{U''_\mu(x_\mu)}\right)^{\frac{3}{2}}(1+o_{\mathcal{I}}(1)).
\end{equation}
{\bf Step 2.3.} Asymptotic behavior of $\tilde{I}_1$. Let us develop the function $U_\mu$ in the neighborhood of $x_\mu$: $U_\mu(t+x_\mu)=U_\mu(x_0)+\frac{1}{2}U''_\mu(x_\mu)t^2+\frac{1}{6}U_\mu^{(3)}(x_\mu)t^3+\frac{1}{24}U_\mu^{(4)}(y_\mu(t))t^4$
where $y_\mu(t)\in[x_\mu-\rho,x_\mu+\rho]$. The regularity of $U_\mu(x)$ with respect to both $x$ and $\mu$ implies the existence of some constant $C>0$ independent of $\mu$ which bounds the forth derivative of $U_\mu$ on the integral support. Therefore we have
\[
f'_m(x_\mu)e^{-C\frac{\rho^4}{\epsilon}}J_\rho\leq e^{\frac{2U_\mu(x_\mu)}{\epsilon}} \tilde{I}_1\leq f'_m(x_\mu)e^{C\frac{\rho^4}{\epsilon}}J_\rho,
\]
with $J_\rho=\int_{-\rho}^{\rho}ze^{-\frac{\mathcal{U}_2}{\epsilon}z^2-\frac{\mathcal{U}_3}{3\epsilon}z^3}dz$ and $\mathcal{U}_k=U_\mu^{(k)}(x_\mu)$. Since $|e^{-x}-1+x-\frac{x^2}{2}|\leq |x|^3e^{|x|}$, we deduce that, for any $z\in[-\rho;\rho]$:
\[
\left|e^{-\frac{\mathcal{U}_3}{3\epsilon}z^3}-1+\frac{\mathcal{U}_3z^3}{3\epsilon}-\frac{\mathcal{U}_3^2z^6}{18\epsilon^2}\right|\leq\left|\frac{\mathcal{U}_3}{3}\right|^3\frac{\rho^9}{\epsilon^3}e^{\frac{\left|\mathcal{U}_3\right|\rho^3}{3\epsilon}}.
\]
We define $m_\rho(l)=\int_{-\rho}^\rho z^l e^{-\frac{\mathcal{U}_2}{\epsilon}z^2}dz$ and $n_\rho(l)=\int_{0}^\rho |z|^l e^{-\frac{\mathcal{U}_2}{\epsilon}z^2}dz$. Some estimation of the integral $J_\rho$ points out directly:
\begin{align*}\left|J_\rho-m_\rho(1)+\frac{\mathcal{U}_3}{3\epsilon}m_\rho(4)-\frac{\mathcal{U}_3^2}{18\epsilon^2}m_\rho(7)\right|
\leq2\left|\frac{\mathcal{U}_3}{3}\right|^3\frac{\rho^9}{\epsilon^3}e^{\frac{\left|\mathcal{U}_3\right|\rho^3}{3\epsilon}}n_\rho(1).
\end{align*}
Symmetry arguments permit easily to deduce that $m_\rho(1)=m_\rho(7)=0$. Finally it suffices to compute $m_\rho(4)$ and $n_\rho(1)$. To this end we introduce the change of variable $u:=\frac{\mathcal{U}_2}{\epsilon}z^2$ and let $\rho^2/\epsilon$ tend to infinity:
\[
m_\rho(4)=\frac{3\sqrt{\pi}}{4}\left(\frac{1}{U''_\mu(x_\mu)}\right)^{\frac{5}{2}}\epsilon^{\frac{5}{2}}(1+o_{\mathcal{I}}(1))\quad\mbox{and}\ n_\rho(1)=\frac{\epsilon}{2 U''_\mu(x_\mu)}(1+o_{\mathcal{I}}(1)).
\]
To sum up: if $\frac{\rho^{18}}{\epsilon^7}\to 0$ (that is $\frac{\rho^9}{\epsilon^2}=o\{\epsilon^{\frac{3}{2}}\}$) then
\begin{equation}\label{step23}	\tilde{I}_1=-\sqrt{\pi}f'_m(x_\mu)\frac{U^{(3)}_\mu(x_\mu)}{4}\left(\frac{1}{U''_\mu(x_\mu)}\right)^{\frac{5}{2}}\epsilon^{\frac{3}{2}}e^{-\frac{2U_\mu(x_\mu)}{\epsilon}}(1+o_{\mathcal{I}}(1)).
\end{equation}
{\bf Step 2.4.} Asymptotic behavior of $\tilde{I}_0$. Let us first study the following integral
\begin{eqnarray*}	I'_0
	&=&\int_{-\rho}^{\rho}\exp\left[-\frac{U''_\mu(x_\mu)}{\epsilon}z^2-\frac{U_\mu^{(3)}(x_\mu)}{3\epsilon}z^3-\frac{U_\mu^{(4)}(x_\mu)}{12\epsilon}z^4\right]dz
	\end{eqnarray*}
We recall the usual notations $\mathcal{U}_k=U^{(k)}_\mu(x_\mu)$. The arguments are similar to those used in Step 2.3. Since $\left|e^{-u}-1+u-\frac{u^2}{2}\right|\leq|u|^3e^{|u|}$, for any $z\in[-\rho;\rho]$ we get
\[
\left|e^{-\frac{\mathcal{U}_3}{3\epsilon}z^3-\frac{\mathcal{U}_4}{12\epsilon}z^4}-1+\frac{\mathcal{U}_3}{3\epsilon}z^3+\frac{\mathcal{U}_4}{12\epsilon}z^4-\frac{1}{2}\left(\frac{\mathcal{U}_3}{3\epsilon}z^3+\frac{\mathcal{U}_4}{12\epsilon}z^4\right)^2\right|\leq C\rho^3.
\]
Adopting the same notations as in Step 2.3 and using symmetry properties, the following bound (uniform with respect to the parameter $\mu$) yields
\[
\left|I'_0-m_\rho(0)+\frac{\mathcal{U}_4}{12\epsilon}m_\rho(4)-\frac{1}{2}\left(\frac{\mathcal{U}_3}{3\epsilon}\right)^2m_\rho(6)-\frac{1}{2}\left(\frac{\mathcal{U}_4}{12\epsilon}\right)^2m_\rho(8)\right|
\leq C\rho^3m_\rho(0).
\]
By the usual change of variable $u:=\frac{U''_\mu(x_\mu)}{\epsilon}z^2$ we emphasize some asymptotic estimation of $I'_0$ as $\rho^2/\epsilon\to\infty$ and $\rho^3/\epsilon\to 0$:
\[
I'_0=\sqrt{\frac{\pi\epsilon}{U''_\mu(x_\mu)}}\left\{1-\frac{U^{(4)}_\mu(x_\mu)}{16U''_\mu(x_\mu)^2}\epsilon+\frac{5U^{(3)}_\mu(x_\mu)^2}{48U''_\mu(x_\mu)^3}\epsilon+o_{\mathcal{I}}(\epsilon)\right\}
\]
We apply the mean value theorem to the function $U_\mu$:
\[
U_\mu(x_\mu+z)=U_\mu(x_\mu)+\frac{\mathcal{U}_2}{2}z^2+\frac{\mathcal{U}_3}{6}z^3+\frac{\mathcal{U}_4}{24}z^4+\frac{1}{120}U^{(5)}(y_\mu(t))z^5,
\]	
with $y_\mu(t)\in[x_\mu-\rho, x_\mu+\rho]$ and $\vert z\vert\le\rho$. From this equality we deduce an estimation of the distance between the integrals $ \tilde{I}_0$ and $I'_0$.\\
We denote by $\mathcal{D}=e^{\frac{2U_\mu(x_\mu)}{\epsilon}}\tilde{I}_0-f_m(x_\mu)I'_0$ this distance. Then there exists some constant $C>0$ independent of $\mu$ and $m$ such that
\begin{align*}
\left|\mathcal{D}\right|&\leq|f_m(x_\mu)|\int_{-\rho}^{\rho}e^{-\frac{\mathcal{U}_2}{\epsilon}z^2-\frac{\mathcal{U}_3}{3\epsilon}z^3-\frac{\mathcal{U}_4}{12\epsilon}z^4}
\left|1-e^{-\frac{1}{60\epsilon}U^{(5)}(y_{\mu}(z+x_\mu)z^5}\right|dz\\
&\leq\frac{\vert f_m(x_\mu)\vert C\lambda}{60\epsilon}\rho^5\int_{-\rho}^{\rho}e^{-\frac{\mathcal{U}_2}{\epsilon}z^2-\frac{\mathcal{U}_3}{3\epsilon}z^3-\frac{\mathcal{U}_4}{12\epsilon}z^4
+\frac{1}{60\epsilon}\left|U^{(5)}(y_\mu(z+x_\mu))z^5\right|}dz.
\end{align*}
If both conditions $\rho^2/\epsilon\to\infty$ and $\rho^3/\epsilon\to 0$ are satisfied then the integral term in the preceding inequality is obviously equivalent to $\sqrt{\frac{\pi\epsilon}{U''_\mu(x_\mu)}}$. The following equivalence holds for the initial integral $\tilde{I}_0$: under the assumption that $\frac{\rho^5}{\sqrt{\epsilon}}=o\left(\epsilon^{\frac{3}{2}}\right)$, we get
$\left|\mathcal{D}\right|=o_{\mathcal{IM}}\left(\epsilon^{\frac{3}{2}}\right)$ and consequently 
\begin{equation}\label{step24}
\tilde{I}_0=e^{-\frac{2U(x_\mu)}{\epsilon}}\sqrt{\frac{\pi\epsilon}{\mathcal{U}_2}}\left\{1-\frac{\mathcal{U}_4}{16\ \mathcal{U}_2^2}\epsilon+\frac{5\ \mathcal{U}_3^2}{48\ \mathcal{U}_2^3}\epsilon+o_{\mathcal{IM}}(\epsilon)\right\}.
\end{equation}
{\bf Step 3.} To sum up: in Step 1, we proved that it suffices to estimate the integral $I_1$ which can be split into 4 terms. Each of them has been estimated in equations \eqref{eq:step21}, \eqref{eq:step22}, \eqref{step23} and \eqref{step24}. The whole integral has the asymptotic equivalence \eqref{eq:lem:inter2-1} as soon as  $\rho^3/\epsilon\to 0$, $\rho^{18}/\epsilon^7\to 0$ and $\rho^5/\epsilon^2\to 0$. The particular choice  $\rho=\epsilon^{\frac{9}{20}}$ fulfills all these conditions.
\end{proof}
We can extend the statement of the preceding lemma to integrals with unbounded supports.
\begin{lem}\label{lem:inter3}
Let $U$ and $G$ be two $\mathcal{C}^\infty(\mathbb{Re})$-continuous functions. We define $U_\mu=U+\mu G$ for the parameter $\mu$ belonging to some compact interval $\mathcal{I}$ of $\mathbb{Re}$ We assume that $U_\mu(t)\geq t^2$ for $|t|$ larger than some $R$ independent of $\mu$ and that $U_\mu$ admits a unique global minimum at $x_\mu$ with  $U''_\mu(x_\mu)>0$. Let $f_m$ be a $\mathcal{C}^3$-continuous function depending on some parameter $m$ which belongs to a compact set $\mathcal{M}$. Furthermore we assume that there exists some constant $\lambda>0$ such that $\left|f_m(t)\right|\leq \exp \left[\lambda| U_\mu(t)|\right]$ for all $t\ge R$, $\mu\in \mathcal{I}$, $m\in\mathcal{M}$ and $\vert f_m^{(i)}\vert$ is locally bounded uniformly with respect to the parameter $m\in \mathcal{M}$ for $0\le i\le 3$.  Then the following asymptotic result holds as $\epsilon$ tends to $0$:
\begin{equation}\label{eq:lem:inter3}
\int_{\mathbb{Re}}f_m(t)e^{\frac{-2U_\mu(t)}{\epsilon}}dt=e^{-\frac{2U_\mu(x_\mu)}{\epsilon}}\sqrt{\frac{\pi\epsilon}{\mathcal{U}_2}}\Big\{f_m(x_\mu)
+\gamma_0(\mu)\epsilon+o_{\mathcal{IM}}(\epsilon)\Big\},
\end{equation}
where $\gamma_0(\mu)$ is defined by \eqref{eq:lem:inter2-2} and $o_{\mathcal{IM}}(\epsilon)/\epsilon$ converges to $0$ as $\epsilon\to 0$ uniformly with respect to the parameters $m$ and $\mu$.
\end{lem}
\begin{proof} Let $R>0$ such that $U_\mu(t)\geq t^2$ for $t\ge R$.
The initial integral can be split into two integrals: the first one denoted by $I_1$ concerns the compact support $[-R,R]$ and the other one $I_2$ concerns the complementary support. For $I_1$ it suffices to apply Lemma \ref{lem:Annexe:inter2} in order to get the asymptotic development. It remains then to prove that $I_2$ is negligible with respect to $I_1$ that is  $I_2=o_{\mathcal{IM}}\left\{\epsilon^{\frac{3}{2}}e^{-\frac{2U_\mu(x_\mu)}{\epsilon}}\right\}$. Using the change of variable $t:=\left(\frac{2}{\epsilon}-\lambda\right)^{-\frac{1}{2}}s$ the following bound holds:
\[
|I_2|\le 2\int_{R}^{+\infty}\exp\left[t^2\left(\lambda-\frac{2}{\epsilon}\right)\right]dt\le 2\sqrt{\frac{\epsilon}{2-\lambda\epsilon}}\int_{R\sqrt{\frac{2-\lambda\epsilon}{\epsilon}}}^{+\infty}\exp\left[-s^2\right]ds.
\]
Lemma \ref{lem:Annexe1} permits to prove as claimed that $I_2$ can be neglected.
\end{proof}
Lemma \ref{lem:inter3} can be applied to particular functions $f_m$. 
\begin{lem}\label{lem:Annexe2}
Let $U$ and $G$ be two $\mathcal{C}^{\infty}(\mathbb{Re})$-continuous functions. We define $U_\mu=U+\mu G$ with $\mu$ belonging to some compact interval $\mathcal{I}$ of $\mathbb{Re}$. We assume that $U_\mu(t)\geq t^2$ for $|t|$ larger than some $R$ independent of $\mu$ and that $U_\mu$ admits a unique global minimum at $x_\mu$ with  $U''_\mu(x_\mu)>0$. Let $f_m$ be a $\mathcal{C}^3$-continuous function depending on some parameter $m$ which belongs to a compact set $\mathcal{M}$. Furthermore we assume that there exists some constant $\lambda>0$ such that $\left|f_m(t)\right|\leq \lambda| U_\mu(t)|$ for all $t\ge R$, $\mu\in \mathcal{I}$, $m\in\mathcal{M}$ and that $\vert f_m^{(i)}\vert$ is locally bounded uniformly with respect to $m\in\mathcal{M}$ for $0\le i \le 3$. Then, for any $n\ge 1$ and asymptotically as $\epsilon\to 0$ we obtain the estimate
\begin{equation*}\label{eq:lem:inter4}
\frac{\int_{\mathbb{R}}t^ne^{f_m(t)}e^{\frac{-2U_\mu(t)}{\epsilon}}dt}{\int_{\mathbb{R}}e^{f_m(t)}e^{\frac{-2U_\mu(t)}{\epsilon}}dt}=x_\mu^n-\frac{nx_\mu^{n-2}}{4\,\mathcal{U}_2}\left[x_\mu\,\frac{\mathcal{U}_3}{\mathcal{U}_2}-n+1-2x_\mu f'_m(x_\mu)\right]\epsilon+o_{\mathcal{IM}}(\epsilon),
\end{equation*}
where $\mathcal{U}_i=U_\mu^{(i)}(x_\mu)$ and $o_{\mathcal{IM}}(\epsilon)/\epsilon$ converges to $0$ as $\epsilon\to 0$ uniformly with respect to the parameters $m$ and $\mu$.
\end{lem}
\begin{proof} We just apply two times Lemma \ref{lem:inter3}: the first time to the denominator $D^\epsilon$ that is for the function $t\to e^{f_m(t)}$ and the second time to the numerator $N^\epsilon$ for the function $t\to t^n e^{f_m(t)}$. The following asymptotic result holds
\begin{equation}\label{ident}
D^\epsilon=e^{-\frac{2U_\mu(x_\mu)}{\epsilon}}\sqrt{\frac{\pi\epsilon}{\mathcal{U}_2}}e^{f_m(x_\mu)}\Big\{1
+\hat{\gamma}_d\epsilon+o_{\mathcal{IM}}(\epsilon)\Big\}
\end{equation}
where
\[
\hat{\gamma}_d=\left(\frac{5\,\mathcal{U}_3^2}{48\,\mathcal{U}_2^3}-\frac{\mathcal{U}_4}{16\,\mathcal{U}_2^2}\right)-f'_m(x_\mu)\frac{\mathcal{U}_3}{4\,\mathcal{U}_2^2}+\Big(f''_m(x_\mu)+f'_m(x_\mu)^2\Big)\frac{1}{4\,\mathcal{U}_2}.
\]
The numerator normalized by $x_\mu^n$ i.e. $N^\epsilon/x_\mu^n$ satisfies some similar identity as $D^\epsilon$, namely \eqref{ident} with $\hat{\gamma}_d$ replaced by $\hat{\gamma}_n$:
\begin{align*}
\hat{\gamma}_n&=\left(\frac{5\,\mathcal{U}_3^2}{48\, \mathcal{U}_2^3}-\frac{\mathcal{U}_4}{16\,\mathcal{U}_2^2}\right)-\left(\frac{n}{x_\mu}+f'_m(x_\mu)\right)\frac{\mathcal{U}_3}{4\,\mathcal{U}_2^2}\\
&+\left(\frac{n(n-1)}{x_\mu^2}+2\frac{n}{x_\mu}f'_m(x_\mu)+f''_m(x_\mu)+f'_m(x_\mu)^2\right)\frac{1}{4\,\mathcal{U}_2}.
\end{align*}
The estimation of the ratio is then a classical exercise of asymptotic analysis.
\end{proof}
The next lemmas are generalizations of Lemma \ref{lem:inter3} and Lemma \ref{lem:Annexe2} to functions $G$ depending on the small parameter $\epsilon$.
\begin{lem}\label{lem:inter-dern}
Let $U$ and $G$ be two $\mathcal{C}^{\infty}(\mathbb{Re})$-continuous functions such that $U(t)\geq t^2$ for $|t|$ large enough and $|G(t)|\le \lambda |U(t)|+C$ for some constants $\lambda>0$ and $C>0$. Moreover we assume that $U$ admits some unique global minimum reached at $x_0$ with $U''(x_0)>0$. 
For any sequence $(\eta_\epsilon)_\epsilon$ satisfying $\lim_{\epsilon\to 0}\eta_\epsilon=0$ and $\lim_{\epsilon\to 0}\epsilon/\eta_\epsilon=0$ we define
$\uem=U+\eta_\epsilon \mu G$ depending on the parameter $\mu$ which belongs to some compact interval $\mathcal{I}$ of $\mathbb{Re}$. Let $f$ a $\mathcal{C}^3$-continuous function such that $|f(t)|\le e^{\lambda|U(t)|}$ for all $|t|$ large enough and such that $\vert f_m^{(i)}\vert$ is locally bounded uniformly with respect to $m\in\mathcal{M}$ for $0\le i\le 3$.
Then, there exists $\epsilon_0>0$ such that the potential $\uem$ admits a unique global minimum reached at $\xem$ for all $\epsilon\le \epsilon_0$. Furthermore the following asymptotic results hold
\begin{equation}
 \label{first}
\xem=x_0-\mu\frac{G'(x_0)}{U''(x_0)}\eta_\epsilon+o_{\mathcal{I}}(\eta_\epsilon)
\end{equation}
\begin{equation}\label{eq:lem:inter-dern}
\int_{\mathbb{R}}f(t)e^{-\frac{2\uem(t)}{\epsilon}}dt=\sqrt{\frac{\pi\epsilon}{U''(x_0)}}e^{-\frac{2\uem(\xem)}{\epsilon}}\Big(f(x_0)+\gamma_\mu\eta_\epsilon+o_{\mathcal{I}}(\eta_\epsilon) \Big),
\end{equation}
where
\[
\gamma_\mu=\frac{\mu}{2U''(x_0)}\Big( -2f'(x_0)G'(x_0)-f(x_0)G''(x_0)+f(x_0)\frac{U^{(3)}(x_0)G'(x_0)}{U''(x_0)} \Big),
\]
and $o_{\mathcal{I}}(\eta_\epsilon)/\eta_\epsilon$ tends to $0$ as $\epsilon\to 0$ uniformly with respect to the parameter $\mu$.

\end{lem}
\begin{proof} Let us first prove that the potential $\uem(x)$ admits a unique minimum for $x=\xem$ with $\lim_{\epsilon\to 0}\xem=x_0$. By the definitions of $(\eta_\epsilon)_\epsilon$ and $\uem$, the following convergence holds
\begin{equation}\label{fin1}
\lim_{\epsilon\to 0}\uem(x_0)=U(x_0).
\end{equation}
Since $x_0$ is the unique global minimum of $U$, for any small $R>0$ there exists $\rho_R>0$ such that $\inf_{x\in[x_0-R,x_0+R]^c}U(x)>U(x_0)+\rho_R$. We deduce the existence of two small constants $\rho'_R$ and $\epsilon_0$ such that  
\begin{align}
 \label{fin2}
\uem(x)\ge (1-\mu\lambda\eta_\epsilon)U(x)-\eta_\epsilon\mu C\ge U(x_0)+\rho'_R,
\end{align}
for all $\epsilon\le\epsilon_0$ and $x\in [x_0-R,x_0+R]^c$. By \eqref{fin1} and \eqref{fin2} we obtain: for any $R>0$ the global minimum of the parametrized potential $\uem$ is reached in the interval $x\in [x_0-R,x_0+R]$ provided that $\epsilon$ is small enough (uniformly with respect to $\mu$). Moreover this global mimimum is unique. Indeed $U''(x_0)>0$ and the regularity of $U$ implies that $U''(x)>0$ for all $x$ in some small neighborhood of $x_0$. Since $\uem''$ converges towards $U''$ as $\epsilon\to 0$ uniformly on each compact subset of $\mathbb{Re}$, we obtain that $\uem''(x)>0$ for all $x\in [x_0-R,x_0+R]$ provided that $R$ and $\epsilon$ are small enough. The minimum is actually unique, we denote its localization $\xem$ and point out that, for $\epsilon$ small, $\uem''(\xem)>0$ uniformly with respect to $\mu$.\\
Let us determine $\xem$. By applying the mean value theorem to $\uem$, we get
\[
0=\uem'(\xem)=U'(x_0)+\mu\eta_\epsilon G'(x_0)+\uem''(\tilde{x})(\xem-x_0),
 \]
where $\tilde{x}$ is in between $x_0$ and $\xem$. Since the second derivative is continuous,   $\uem''(\tilde{x})$ is uniformly bounded. Moreover $U'(x_0)=0$. Consequently $\xem-x_0=\mathcal{O}_{\mathcal{I}}(\eta_\epsilon)$. Using the same argument for the second order asymptotic development of $\uem'(\xem)$, that is
\[
 0=U'(x_0)+\mu\eta_\epsilon G'(x_0)+\Big( U''(x_0)+\mu\eta_\epsilon G''(x_0) \Big) (\xem-x_0)+\frac{\uem^{(3)}(\tilde{x})}{2}(\xem-x_0)^2,
\]
we obtain the announced estimate \eqref{first}. Finally let us prove the estimate \eqref{eq:lem:inter-dern}. The statement of Lemma \ref{lem:inter3} can be applied to $\uem$ since the asymptotic result \eqref{eq:lem:inter3} is uniform with respect to the parameter $\mu$. So it suffices to consider the case when $\mu$ is replaced by $\mu\eta_\epsilon$. We immediately obtain
\begin{equation}\label{last}
 \int_\mathbb{Re}f(t)e^{-\frac{\uem(t)}{\epsilon}}dt=\sqrt{\frac{\pi\epsilon}{\uem''(\xem)}}f(\xem)e^{-\frac{\uem(\xem)}{\epsilon}}\Big(1+o_{\mathcal{I}}(\eta_\epsilon)\Big).
\end{equation}
It remains to approximate $ f(\xem)$ and $\uem''(\xem)$ using \eqref{first}. Due to the regularity of both $f$ and $U$, the following developments hold  
 \[
  f(\xem)=f(x_0)-\mu\eta_\epsilon f'(x_0)\frac{G'(x_0)}{U''(x_0)}+o_{\mathcal{I}}(\eta_\epsilon),
 \]
\[
 \uem''(\xem)=U''(x_0)+\mu\eta_\epsilon\Big(G''(x_0)-U^{(3)}(x_0)\frac{G'(x_0)}{U''(x_0)}\Big)+o_{\mathcal{I}}(\eta_\epsilon).
\]
The statement of Lemma \ref{lem:inter-dern} is obtained just by combination of the two preceding asymptotics and \eqref{last}.
\end{proof}
We are now able to present a statement similar to Lemma \ref{lem:Annexe2} for some potential $U_\mu$ depending on the small parameter $\epsilon$. It suffices to consider a ratio of two integral terms. Then an immediate application of Lemma \ref{lem:inter-dern} leads to the following result.
\begin{lem}\label{lem:Annexe3}
Let $U$ and $G$ be two $\mathcal{C}^{\infty}(\mathbb{Re})$-continuous functions such that $U(t)\geq t^2$ for $|t|$ large enough and $|G(t)|\le \lambda |U(t)|+C$ for some constants $\lambda>0$ and $C>0$. Moreover we assume that $U$ admits some unique global minimum reached at $x_0$ with $U''(x_0)>0$. 
For any sequence $(\eta_\epsilon)_\epsilon$ satisfying $\lim_{\epsilon\to 0}\eta_\epsilon=0$ and $\lim_{\epsilon\to 0}\epsilon/\eta_\epsilon=0$ we define
$\uem=U+\eta_\epsilon \mu G$ depending on the parameter $\mu$ which belongs to some compact interval $\mathcal{I}$ of $\mathbb{Re}$. Let $f$ a $\mathcal{C}^3$-continuous function such that $|f(t)|\le e^{\lambda|U(t)|}$ for all $|t|$ large enough and such that $\vert f_m^{(i)}\vert$ is locally bounded uniformly with respect to $m\in\mathcal{M}$ for $0\le i\le 3$. Then as $\epsilon\to 0$, we obtain the following estimate
\begin{equation}\label{eq:lem:A3}
\frac{\int_{\mathbb{Re}}f(t)e^{-\frac{2\uem(t)}{\epsilon}}dt}{\int_{\mathbb{Re}}e^{-\frac{2\uem(t)}{\epsilon}}dt}=f(x_0)-\frac{f'(x_0)G'(x_0)}{U''(x_0)}\eta_\epsilon+o_{\mathcal{I}}(\eta_\epsilon)
\end{equation}
where $o_{\mathcal{I}}(\eta_\epsilon)/\eta_\epsilon$ tends to $0$ as $\epsilon\to 0$ uniformly with respect to the parameter $\mu$.
\end{lem}
\begin{rem}
 \label{lem:rem1}
The statements of Lemmas \ref{lem:Annexe:inter1}-\ref{lem:Annexe3} can be easely generalized, replacing the parametrized function $U_\mu=U+\mu G$ by $U_\mu=U+\sum_{i=1}^k \mu_i G_i$ where $\mu=(\mu_1,\ldots,\mu_k)\in \mathcal{I}_1\times\ldots\times \mathcal{I}_k$. The convergence results are then uniform with respect to all parameters.
\end{rem}

    \begin{small}
 \bibliographystyle{plain}
 \bibliography{tugautherrmann.bbl}
\end{small}
\end{document}